\documentclass{IEEEtran}
\usepackage{amsmath}
\usepackage{amssymb}
\usepackage{amsthm}
\usepackage[pdftex]{graphicx}
\usepackage{geometry}
\usepackage{epsfig}
\usepackage{epstopdf}
\usepackage{mathtools}
\usepackage{relsize}
\usepackage{hyperref}
\hypersetup{
    colorlinks=true,     
    linkcolor=blue,     
    citecolor=blue,     
}
\usepackage[numbers,sort&compress]{natbib}
\usepackage{enumitem}
\usepackage{algorithm,algpseudocode}
\usepackage{fullpage}
\usepackage{setspace}
\usepackage{caption}
\usepackage{subcaption} 
\usepackage{xspace}
\usepackage{url}
\usepackage{multirow}
\usepackage{color}
\usepackage{authblk}

\newtheorem{theorem}{Theorem}
\newtheorem{lemma}[theorem]{Lemma}
\newtheorem{assumption}[theorem]{Assumption}
\newtheorem{remark}[theorem]{Remark}
\newtheorem{definition}[theorem]{Definition}
\newtheorem{corollary}[theorem]{Corollary}
\newtheorem{proposition}[theorem]{Proposition}
\newtheorem{example}[theorem]{Example}
\newtheorem{conjecture}[theorem]{Conjecture}

\let\oldtheorem\theorem
\renewcommand{\theorem}{\oldtheorem\normalfont}

\let\oldlemma\lemma
\renewcommand{\lemma}{\oldlemma\normalfont}

\let\oldassumption\assumption
\renewcommand{\assumption}{\oldassumption\normalfont}

\let\oldremark\remark
\renewcommand{\remark}{\oldremark\normalfont}

\let\olddefinition\definition
\renewcommand{\definition}{\olddefinition\normalfont}

\let\oldcorollary\corollary
\renewcommand{\corollary}{\oldcorollary\normalfont}

\let\oldproposition\proposition
\renewcommand{\proposition}{\oldproposition\normalfont}

\let\oldexample\example
\renewcommand{\example}{\oldexample\normalfont}

\let\oldconjecture\conjecture
\renewcommand{\conjecture}{\oldconjecture\normalfont}

\newcommand{\Set}[2]{\lbrace #1 : #2 \rbrace}

\newcommand{\uset}[2]{\ensuremath{\underset{#1}{#2}}}

\DeclarePairedDelimiter\abs{\lvert}{\rvert}%

\newcommand{\prob}[1]{\mathbb{P}\left\lbrace{#1}\right\rbrace}
\newcommand{\expect}[1]{\mathbb{E}\left[{#1}\right]}
\newcommand{\expectation}[2]{\mathbb{E}_{#1}\hspace{-0.03in}\left[{#2}\right]}
\newcommand{\stdev}[1]{\sigma\left[{#1}\right]}

\newcommand{\D}{\mathcal{D}}
\newcommand{\E}{\mathcal{E}}
\newcommand{\G}{\mathcal{G}}

\newcommand{\Gres}{\mathcal{G}^{\text{res}}}
\newcommand{\Rg}{\mathcal{R}}
\newcommand{\Resa}{\mathcal{R}es(\alpha)}
\newcommand{\V}{\mathcal{V}}
\newcommand{\W}{\mathcal{W}}

\newcommand{\random}{\omega}

\newcommand{\dtiw}{\tilde{d}_i(\omega)}
\newcommand{\dtjw}{\tilde{d}_j(\omega)}

\newcommand{\fmax}{f^{\textup{max}}}
\newcommand{\po}{p^{0}}

\newcommand{\pbojw}{\bar{p}^0_j(\random)}
\newcommand{\pw}{p(\random)}
\newcommand{\piw}{p_i(\random)}
\newcommand{\pt}{p^{\textup{T}}}
\newcommand{\ptw}{p^{\textup{T}}(\random)}
\newcommand{\ptiw}{p^{\textup{T}}_i(\random)}
\newcommand{\pmin}{p^{\textup{min}}}
\newcommand{\pmax}{p^{\textup{max}}}
\newcommand{\pminiw}{p^{\textup{min}}_i(\random)}
\newcommand{\pmaxiw}{p^{\textup{max}}_i(\random)}
\newcommand{\pminjw}{p^{\textup{min}}_j(\random)}
\newcommand{\pmaxjw}{p^{\textup{max}}_j(\random)}
\newcommand{\pol}{p^{0,\textup{L}}}
\newcommand{\pou}{p^{0,\textup{U}}}
\newcommand{\rp}{r^+}
\newcommand{\rpu}{r^{+,\textup{max}}}
\renewcommand{\rm}{r^-}
\newcommand{\rmu}{r^{-,\textup{max}}}
\newcommand{\sw}{s(\random)}
\newcommand{\swk}{s(\random_k)}
\newcommand{\slw}{s^{\textup{L}}(\random)}
\newcommand{\suw}{s^{\textup{U}}(\random)}
\newcommand{\sigmadw}{\Sigma_d(\random)}
\newcommand{\tw}{\theta(\random)}
\newcommand{\tiw}{\theta_i(\random)}
\newcommand{\tjw}{\theta_j(\random)}

\newcommand{\xl}{x^{\textup{L}}}
\newcommand{\xu}{x^{\textup{U}}}

\newcommand{\casenum}[1]{\texttt{case{#1}}}

\makeatletter
\newcommand*\bigcdot{\mathpalette\bigcdot@{.6}}
\newcommand*\bigcdot@[2]{\mathbin{\vcenter{\hbox{\scalebox{#2}{$\m@th#1\bullet$}}}}}
\makeatother

\makeatletter
\newcommand{\algmargin}{\the\ALG@thistlm}
\makeatother
\newlength{\whilewidth}
\algdef{SE}[parWHILE]{parWhile}{EndparWhile}[1]
  {\parbox[t]{\dimexpr\linewidth-\algmargin}{%
     \hangindent\whilewidth\strut\algorithmicwhile\ #1\ \algorithmicdo\strut}}{\algorithmicend\ \algorithmicwhile}%
\algnewcommand{\parState}[1]{\State%
  \parbox[t]{\dimexpr\linewidth-\algmargin}{\strut #1\strut}}

\title{Stochastic DC Optimal Power Flow \\ With Reserve Saturation\thanks{This research is supported by the Department of Energy, Office of Science, Office of Advanced Scientific Computing Research, Applied Mathematics program under Contract Number DE-AC02-06CH11347.}}

\author[1]{Rohit Kannan}
\author[2]{James R. Luedtke}
\author[3]{Line A. Roald}
\affil[1]{Wisconsin Institute for Discovery, University of Wisconsin-Madison, Madison, WI, USA. \protect\\ E-mail: rohit.kannan@wisc.edu}
\affil[2]{Department of Industrial \& Systems Engineering and Wisconsin Institute for Discovery, University of Wisconsin-Madison, Madison, WI, USA. E-mail: jim.luedtke@wisc.edu}
\affil[3]{Department of Electrical and Computer Engineering, University of Wisconsin-Madison, Madison, WI, USA. \protect\\ Corresponding Author. E-mail: roald@wisc.edu}

\begin{document}


\twocolumn[
  \begin{@twocolumnfalse}
    \maketitle
    \begin{abstract}
We propose an optimization framework for stochastic optimal power flow with uncertain loads and renewable generator capacity. Our model follows previous work in assuming that generator outputs respond to load imbalances according to an affine control policy, but introduces a model of {\it saturation} of generator reserves by assuming that when a generator's target level hits its limit, it abandons the affine policy and produces at that limit. 
This is a particularly interesting feature in models where wind power plants, which have uncertain upper generation limits, are scheduled to provide reserves to balance load fluctuations.
The resulting model is a nonsmooth nonconvex two-stage stochastic program, and we use a stochastic approximation method to find stationary solutions to a smooth approximation. 
Computational results on 6-bus and 118-bus test instances demonstrates that by considering the effects of saturation, our model can yield solutions with lower expected generation costs (at the same target line violation probability level) than those obtained from a model that enforces the affine policy to stay within generator limits with high probability. 
    \end{abstract}

  \end{@twocolumnfalse}
]

\section{Introduction}
Large shares of renewable energy increases the variability and uncertainty in power grid operations, and frequently lead to a larger demand for balancing energy through generation reserves. 
Understanding and counteracting potentially adverse effects of this uncertainty requires models that accurately capture its impact on the network.
Ignoring the effect of uncertainties while making dispatching decisions can result in unsafe operations~\citep{bienstock2014chance}, whereas considering them can significantly improve system security while simultaneously enabling economic efficiency~\citep{roald2016optimal}.
Many approaches to stochastic optimal power flow (OPF) problems have typically relied on affine generation control policies to balance fluctuating power demands, mimicking the actions of the automatic generation control~\citep{vrakopoulou2013,bienstock2014chance,roald2015optimal,roald2016optimal,zhang2016distributionally}. 
These policies require traditional generators to provide a determined fraction of the necessary reserves.
The feasibility of the affine control policy is typically enforced using conservative chance-constrained approximations~\citep{vrakopoulou2013,bienstock2014chance,zhang2016distributionally}, robust constraints \cite{warrington2013}, or by constraining the expected exceedance of determined reserves~\citep{roald2015optimal,roald2016optimal}. 
A key limitation of the affine control policy is that it does not adequately model the behavior of the generators as they reach their upper or lower generation limits~\citep{roald2015optimal,roald2016optimal}.
When the system faces large demand fluctuations, some generators are likely to hit their limits if the affine policy is used, in which case a realistic generator will simply stop providing reserves and maintain a fixed power output. 
Failing to model this behavior may result in conservative results with economically inferior dispatching decisions,
because requiring feasibility of the affine policy forces each generator to maintain too large of a reserve capacity.
This drawback becomes more pronounced when considering reserves from uncertain resources, such as reserves provided by renewable generators themselves \citep{roald2016optimal}, or demand response resources \citep{zhang2016distributionally}. 

To address this drawback, we introduce an new optimization model 
that includes a more realistic and flexible representation of reserve activation and captures the impact of upper and lower generator limits, which we call \emph{reserve saturation}. While this reserve saturation model is an accurate reflection of current system operations, it has, to the best of our knowledge, never before been considered in the context of stochastic optimal power flow. Related work resets the affine control policy through activation of manual reserves \citep{roald2015optimal}, imposes hard limitations on wind power generation~\citep{roald2016optimal}, or use multi-parametric programming as a preprocessing step \citep{vrakopoulou2017optimal}. However, all these methods require the user to pre-specify important aspects of the piecewise affine policies, leading to potentially sub-optimal solutions.
In contrast, we introduce a two-stage stochastic formulation for the DC OPF problem that includes the reserve saturation model in the second-stage, 
which inherently incorporates and enforces power generation limits in conventional and wind generators. By explicitly enforcing the piecewise control policy as an second-stage constraint, the optimization problem is able to identify the optimal generation and reserve allocation considering this behavior, without any pre-specified (and potentially sub-optimal) input.


After introducing the model, we investigate conditions under which it is feasible and derive a stochastic approximation method for solving a smooth version of it to local optimality. 
Finally, we demonstrate the practical benefits of our modeling framework on case studies based on a small 6-bus system and the IEEE 118-bus system.
In particular, we assess the economic and environmental impact of allowing wind power plants to provide significant reserves, and demonstrate empirically that our approach finds solutions that satisfy physical limits with high probability through out-of-sample testing\footnote{For large enough line penalties (to be defined in Sect.~\ref{sec:model}), our solutions are theoretically guaranteed, by Markov's inequality, to satisfy line limits with high probability}.

This paper is organized as follows.
Section~\ref{sec:model} outlines our reserve saturation model within a two-stage stochastic programming framework, and Section~\ref{sec:solnapproach} presents a projected stochastic gradient method for solving an approximation. 
Section~\ref{sec:alternatives} briefly discusses alternative modeling approaches for determining candidate first-stage solutions.
Computational results are reported in Section~\ref{sec:computexp}, and we conclude in Section~\ref{sec:concl}.

\textit{Notation}. We denote vectors by lower case letters and their components using subscripts. 
We let $\text{int}(S)$ denote the interior of a set~$S$, write $(\cdot)_+$ and $(\cdot)_{-}$ to denote $\max\{\cdot,0\}$ and $\min\{\cdot,0\}$, write $\text{logspace}(a,b,n)$ to denote a vector of $n$ logarithmically-spaced points between $10^{a}$ and $10^{b}$ (both inclusive),
and write $\expect{\cdot}$ and $\stdev{\cdot}$ to denote expectation and standard deviation operators.
We do not make a notational distinction between random variables and their realizations.

\section{Optimal power flow with reserve saturation}
\label{sec:model}


We introduce a two-stage stochastic programming model for determining power generation and reserve levels in a power system facing random loads and wind generation uncertainty. 
In the first stage,
the nominal generation levels, reserve capacities and reserve participation factors for each generator are determined. These decisions are taken in advance of observing the random demand and wind generation capacity. 
The second stage models the system response to the observed load and wind generation. This response, which requires the generators to activate reserves to balance the system, is determined by the reserve participation factors from the first-stage of the optimization model. 
In our model, the random loads are uncertain and non-dispatchable,  representing a combination of standard load and non-dispatchable renewable generation. 
We assume that wind power plants are fully dispatchable, except that their output is capped by the random available capacity.

A novel feature of our model is that we explicitly model {\it generator saturation} in the second-stage formulation, which occurs when the output of a generator, as determined by its nominal generation level, participation factor, load imbalance, and the control policy reaches its upper or lower generation limit. A generator that reaches its lower/upper limit continues to produce at that limit, and any additional balancing energy must be provided by the remaining generators that have not yet reached saturation. In this model, generators are allowed to exceed their scheduled reserve capacity, but we assume the system operator pays a higher price for doing so. 
While generation limits are satisfied by virtue of our modeling framework, we use a penalty on the the expected violation of line limits to obtain a solution that satisfies the line limits with high probability.
The objective in our model is hence to minimize the expected generation costs while keeping the expected violation of the line limits small. 



\subsection{Network representation}

We model the network as an undirected connected graph $G = (\V,\E)$, where $\V$  denotes the set of nodes/buses and $\E$ denotes the set of edges/transmission lines. 
The set of wind generators, regular generators, and loads are denoted by $\W$, $\Rg$, and $\D$, respectively, and $\G := \Rg \cup \W$ denotes the set of all generators.
The demand at node $i \in \D$ is a random variable having expected value $d_i$.
We assume for notational convenience that each node in the network houses a load \textit{and} either a wind generator, or a regular generator. It is straightforward to extend the model to include nodes with no/multiple generators/loads.

\subsection{First-stage decisions and constraints}

The first-stage decisions include, for each generator $i \in \G$, the nominal generation levels $\po_i$, the scheduled up- and down-reserve levels $\rp_i$ and $\rm_i$, 
and the generator participation factors for reserves, $\alpha_i$. Note that the set of dispatchable generators includes the wind generators. 
To ensure consistency of our control policy (cf.~\citep{bienstock2014chance,roald2016optimal}), we require that the nominal generation levels satisfy a power balance constraint for the expected value of the demands and lie within pre-specified bounds ($\pol$ and $\pou$):
\begin{align}
\label{eqn:init_bal}
&\sum_{i \in \G} \po_i = \sum_{j \in \D} d_j, \quad \pol \leq \po \leq \pou.
\end{align}
For regular generators, $\pol=\pmin$ and $\pou=\pmax$ represent the upper and lower generation limits. For wind power plants, $\pol$ and $\pou$ represent the maximum and minimum generation that the operator is willing to schedule from that plant.

The up- and down-reserve levels $\rp$ and $\rm$ are constrained to lie within pre-specified limits ($\rpu,\rmu$), and comply with capacity limits for regular generators:
\begin{align}
&\quad 0 \leq \rp \leq \rpu, \quad 0 \leq \rm \leq \rmu, \label{eqn:reserve_bounds} \\
&\po_i + \rp_i \leq \pmax_i, \quad \po_i - \rm_i \geq \pmin_i, \quad \forall i \in \Rg. \label{eqn:reserve_limits}
\end{align}
We assume that the reserve activation follows through the automatic generation control (AGC), where the contribution of each generator is determined through a participation factor \cite{wood2012}. 
The participation factors~$\alpha$ are required to sum to one, and
a subset of generators $\Gres \subset \G$ are required to provide reserves with a minimum participation factor~$\varepsilon$:
\begin{align}
\label{eqn:participation}
&\alpha \geq 0, \quad \sum_{i \in \G}\alpha_i = 1, \quad \alpha_i \geq \varepsilon, \:\: \forall i \in \Gres.
\end{align}
Let $\Resa := \Set{i \in \G}{\alpha_i > 0}$ denote the set of generators with positive participation factors for any choice of~$\alpha$ satisfying constraint~\eqref{eqn:participation}, and note that $\Gres \subset \Resa$.

\subsection{Uncertain parameters and the recourse policy}
We let $\random$ denote the underlying random variables, and assume that we can generate iid samples from its probability distribution.
Let $\dtiw$ represent the random fluctuations in the power demands for $i \in \D$, and $\pminiw$ and $\pmaxiw$, $i \in \W$, represent the minimum and maximum wind generator power outputs. For notational simplicity, we also define $\pmaxiw \equiv \pmax_i$ and $\pminiw \equiv \pmin_i$ for $i \in \Rg$.
Let $\Sigma_d(\random) := \sum_{i \in \D} \dtiw$ denote the net demand fluctuation.

A common assumption in power system modeling~\citep{bienstock2014chance,zhang2016distributionally,vrakopoulou2017optimal, roald2016optimal,wood2012} is that reserve activation to balance the load fluctuation $\Sigma_d(\random)$ follows an affine control policy. This policy adjusts the generation of the regular and wind generators as 
\begin{align}
\label{eqn:affine_policy}
\piw &= \po_i + \alpha_i \Sigma_d(\random), \quad \forall i \in \G,
\end{align}
where $\piw$ denotes the power output of generator $i \in \G$ for a realization of the random variables~$\random$.
While the affine policy satisfies the total power balance constraint by virtue of Eqns.~\eqref{eqn:init_bal} and~\eqref{eqn:participation}, the generation levels~$\piw$ determined by this policy could exceed the generation limits $\pminiw$ and $\pmaxiw$ if the magnitude of the deviation $\Sigma_d(\random)$ is large.
To avoid such violations, existing approaches~\citep{bienstock2014chance,zhang2016distributionally,vrakopoulou2017optimal,
roald2016optimal} impose tight constraints on the probability or expected violation of generation limits, leading to conservative nominal generation levels $\po$ and allocation of the participation factors~$\alpha$.

We propose a more realistic and physically accurate model includes \emph{reserve saturation}. The reserve saturation model allows generators to provide reserves with a determined participation factor \textit{only} until they hit their generation limits, after which other non-saturated generators are required to contribute additional reserves according to their relative participation factors.
To represent this model, we first define \emph{target generation levels} $\ptiw$ as follows: 
\begin{align}
\label{eqn:targets}
\ptiw &= \po_i + \alpha_i \sigmadw + \alpha_i \sw, \quad \forall i \in \G, 
\end{align}
where $\sw$ is a slack reserves variable 
which covers the imbalance incurred by generators that have reached their bounds and are no longer contributing reserves.
For each generator $i \in \G$, the \emph{actual generation levels} $\piw$ are determined using the piecewise-affine policy
\begin{align}
\label{eqn:saturation}
\piw &= \begin{cases}
    \pminiw, & \text{if } \ptiw < \pminiw \\
    \ptiw, & \text{if } \pminiw \leq \ptiw \leq \pmaxiw \\
    \pmaxiw, & \text{if } \ptiw > \pmaxiw .
  \end{cases} 
\end{align}

While the first stage only requires the nominal generation and load to be balanced, the second stage includes the DC power flow constraints for each node~$i \in \V$:
\begin{align}
\label{eqn:flowbal}
\sum_{j \: : \: (i,j) \in \E} \beta_{ij} \left[ \tiw - \tjw \right] &= \piw - d_i - \dtiw,
\end{align} 
where $\tiw$ denotes the phase angle at bus $i \in \V$ and $\beta_{ij}$ ($=\beta_{ji}$) denotes the susceptance in the line $(i,j)$.
Summing Eqn.~\eqref{eqn:flowbal} yields the following total power balance constraint:
\begin{align}
\label{eqn:final_bal}
\sum_{i \in \G} \piw &= \sum_{j \in \D} \left( d_j + \dtjw \right).
\end{align}
For any given values of the generator levels $\piw$, there is a one-dimensional affine space of solutions $\tiw$ to Eqn.~\eqref{eqn:flowbal}.
We assume without loss of generality that the first node is chosen as the reference bus with $\theta_1(\random) \equiv 0$, which, along with Eqn.~\eqref{eqn:flowbal}, implies that there is a unique solution to the phase angles $\tiw$, $i \in \V$, e.g., see Lemma~1.1 of~\citep{bienstock2014chance}.
Line flows $\left[\beta_{ij}(\tiw - \tjw)\right]$ are encouraged to obey line limits by using penalty terms in the objective function.

\vspace*{-0.05in}

\subsection{Solution to the second-stage problem}

We now characterize conditions under which the system of equations~\eqref{eqn:targets},~\eqref{eqn:saturation}, and~\eqref{eqn:final_bal} has a (unique) solution for the second-stage variables $\left(\ptiw, \piw, \sw \right)$ given fixed values for the first-stage variables $(\po,\rp,\rm,\alpha)$.

\begin{theorem}
\label{thm:consis_uniq}
For each value of the first-stage variables satisfying Eqns.~\eqref{eqn:init_bal} to~\eqref{eqn:participation}, the system of equations~\eqref{eqn:targets},~\eqref{eqn:saturation}, and~\eqref{eqn:final_bal} is feasible exactly when $\sigmadw \in D_F(\po,\alpha,\random)$, where $D_F(\po,\alpha,\random) :=$
{
\small
\begin{align*}
&\Bigl[ \sum_{i \in \Resa} \pminiw + \sum_{j \not\in \Resa} \pbojw - \sum_{k \in \D} d_k,  \\
&\quad\quad  \sum_{i \in \Resa} \pmaxiw + \sum_{j \not\in \Resa} \pbojw - \sum_{k \in \D} d_k \Bigr],
\end{align*}
}
$j \not\in \Resa$ is shorthand for $j \in \G \backslash \Resa$, and
{
\small
\begin{align*}
\pbojw &= \textup{median}\left(\po_j, \pminjw, \pmaxjw \right), \quad \forall j \not\in \Resa.
\end{align*}
}
Furthermore, the solution for the $\piw$ variables is always unique, whereas the solution for $\left(\ptiw, \sw \right)$ is unique iff $\sigmadw \in \textup{int}(D_F(\po,\alpha,\random))$.
\end{theorem}

The proof for Theorem \ref{thm:consis_uniq} can be found in Appendix \ref{app:proof}. We henceforth assume that $\sigmadw \in \text{int}(D_F(\po,\alpha,\random))$ for a.e.~realization of~$\omega$ for each value of $\po$ and $\alpha$ satisfying Eqns.~\eqref{eqn:init_bal} to~\eqref{eqn:participation}.
Theorem~\ref{thm:consis_uniq} and its proof then implies that given first-stage decisions $\po$ and $\alpha$ and a realization of the random variables~$\random$, computing \textit{the} recourse solution reduces to a one-dimensional search for the slack reserves~$\sw$.

\subsection{Two-stage stochastic programming model}

We propose the following two-stage stochastic DC-OPF model with reserve saturation:
\begin{alignat}{2}
\label{eqn:stochopf}
&\uset{\rm, \alpha}{\min_{\po,\rp,}} && \displaystyle\sum_{i \in \G} \left[ f_{1,i}(\po_i) + f_{2,i}(\rp_i) + f_{3,i}(\rm_i) \right] + \nonumber \\
& && \quad\quad\quad Q(\po, \rp, \rm, \alpha) \nonumber \\
&\quad \text{s.t.} && \:\: \text{Constraints}~\eqref{eqn:init_bal} \text{ to}~\eqref{eqn:participation}, \tag{P}
\end{alignat}
where $Q(\po,\rp,\rm,\alpha) = \expectation{\random}{q(\po,\rp,\rm,\alpha,\random)}$ denotes the expected second-stage costs with $q(\po,\rp,\rm,\alpha,\random) :=$
{
\begin{alignat}{2}
\label{eqn:recourse}
\uset{\sw,\tw}{\min_{\pw,\ptw,}} \nonumber &\sum_{i \in \G} \big[ && q_{1,i}\left(\piw - \po_i\right) + \big. \\
& && \big. q_{2,i}\left(\piw - (\po_i + \rp_i)  \right) + \big.  \nonumber\\
& && \big. q_{3,i}\left(\piw - (\po_i - \rm_i)  \right) \big] \: + \nonumber \\
&\sum_{(i,j) \in \E} && q_{4,ij}\left(\beta_{ij} \left[ \tiw - \tjw \right]\right) \tag{R} \\
\text{s.t.} \quad & && \hspace*{-0.4in}\text{Constraints}~\eqref{eqn:targets} \text{ to}~\eqref{eqn:flowbal}. \nonumber
\end{alignat}
}
The functions~$f_1$, $f_2$, and $f_3$ in the first-stage objective quantify the cost of nominal power generation and the cost of up- and down-reserve capacities, respectively.
In the second-stage problem, the term involving the function~$q_1$ quantifies the cost of deviating from the generation level decided in the first-stage, representing e.g., a mileage payment to generators, whereas the terms involving the functions~$q_2$ and~$q_3$ correspond to the penalties for using up- and down-reserves beyond the scheduled reserve limits.
Finally, the function~$q_4$ penalizes `large line flows', with the penalty coefficient chosen to trade-off between the cost of power generation and the line flow violation probability. For simplicity, we use a linear weighting approach for balancing the expected generation cost and the expected violation cost. Alternatively, a constraint on the expected violation penalty could be imposed. 
We assume that functions~$f_1$ to~$f_3$ and $q_1$ to~$q_4$ are continuously differentiable with Lipschitz continuous gradients.


\section{Solution approach}
\label{sec:solnapproach}
Modeling reserve saturation introduces bilinear terms (in the expression for the target generation levels~\eqref{eqn:targets}) and nonsmooth nonconvex functions (in the saturation policy~\eqref{eqn:saturation}).
Thus, Problem~\eqref{eqn:stochopf} is a nonsmooth nonconvex two-stage stochastic program, which is in general a challenging problem class to solve even to local optimality. 
However, Theorem~\ref{thm:consis_uniq} implies that for given first-stage decisions, the unique second-stage solution can essentially be computed by solving a one-dimensional equation.
In this section, we further show that by replacing the saturation Eqn.~\eqref{eqn:saturation} with a suitable smooth approximation, partial derivatives of the recourse solution with respect to the first-stage decisions can be computed by solving a linear system.
Therefore, Theorem~\ref{thm:consis_uniq} suggests a sampling-based decomposition approach (i.e., an approach that works in the space of the first-stage variables) for solving an approximation of Problem~\eqref{eqn:stochopf} to obtain stationary solutions.
The remainder of this section proposes a smooth approximation of Problem~\eqref{eqn:stochopf} and a stochastic approximation-based~\citep{ghadimi2016mini,davis2018stochastic} decomposition approach for solving it.


\subsection{Smooth approximation of Problem~\eqref{eqn:stochopf}}

We propose a smooth approximation of Problem~\eqref{eqn:stochopf} to obtain a formulation in which all functions are continuously differentiable.
The nonsmooth saturation function in 
\eqref{eqn:saturation} is approximated by the continuously differentiable function
%
\begin{equation}
\label{eqn:smooth_sat}
\piw = g_{\tau_{sat}}(\ptiw;\pminiw,\pmaxiw)
\end{equation}
where
{
\footnotesize
\begin{align*}
&g_{\tau}(x;\xl,\xu) :=  \\
&\begin{cases}
    \xl, &\text{if } x < \xl - \tau \\
    \xl + \left(x - (\xl - \tau)\right)^2/4\tau, &\text{if } \xl - \tau \leq x \leq \xl + \tau \\
    x, &\text{if } \xl + \tau < x < \xu - \tau \\
    \xu - \left(x - (\xu + \tau)\right)^2/4\tau, &\text{if } \xu - \tau \leq x \leq \xu + \tau \\
    \xu, &\text{if } x > \xu + \tau
  \end{cases},
\end{align*}
}%
and $\tau_{sat} > 0$ is a parameter that controls the approximation quality.
We call the approximation of Problem~\eqref{eqn:stochopf} resulting from this modification `\textit{the smooth approximation}'.

\subsection{Solving the recourse problem of the smooth approximation}

The analysis of Theorem~\ref{thm:consis_uniq} carries over to the smooth approximation because the function~$g_{\tau}$ is monotonically nondecreasing.
Therefore, given values of the first-stage variables and a realization of~$\random$, the \textit{unique} recourse solution can be computed by solving the one-dimensional equation for the slack reserves~$\sw$ that results from substituting Eqns.~\eqref{eqn:targets} and~\eqref{eqn:smooth_sat} into Eqn.~\eqref{eqn:final_bal}.
We solve this equation by bisection.

Denote the (target) power generation levels obtained from Eqns.~\eqref{eqn:targets} and~\eqref{eqn:smooth_sat} with $\sw = 0$ by $\hat{p}^{\text{T}}_i(\random): = \po_i + \alpha_i \sigmadw$ and $\hat{p}_i(\random) := g_{\tau_{sat}}(\hat{p}^{\text{T}}_i(\random);\pminiw,\pmaxiw)$.
Let the residual of Eqn.~\eqref{eqn:final_bal} at these generation levels be denoted by $\delta d(\random) := \sum_{i \in \G} \hat{p}_i(\random) - \sum_{j \in \D} \left( d_j + \dtjw \right)$. If $\delta d(\random)$ is negative, we need to provide up-reserves, whereas if $\delta d(\random)$ is positive, we need to provide down-reserves to balance the overall load for the smooth approximation.
Depending on whether $\delta d(\random) < 0$, or $\delta d(\random) > 0$, we use either $\slw = -\delta d(\random)$ and
{
\small
\begin{align*}
\suw &= \uset{i \in \Resa}{\max} \left\lbrace\alpha_i^{-1}(\pmaxiw + \tau_{sat} - \po_i) - \sigmadw\right\rbrace,
\end{align*}
}
or $\suw = -\delta d(\random)$ and
{
\small
\begin{align*}
\slw &= \uset{i \in \Resa}{\min} \left\lbrace\alpha_i^{-1}(\pminiw - \tau_{sat} - \po_i) - \sigmadw\right\rbrace
\end{align*}
}%
as lower and upper bounds ($s^{\text{L}}$, $s^{\text{U}}$) for the bisection procedure.

\subsection{Computing stochastic gradients for the approximation}

We describe how stochastic gradients of the objective function of Problem~\eqref{eqn:stochopf} are estimated given values of the first-stage variables and a realization of~$\random$.
Given partial derivatives of the (unique) recourse solution with respect to the first-stage decisions, we can compute stochastic gradients of the objective function of the approximation using the chain rule under mild conditions (see Theorem~7.44 of~\citep{shapiro2009lectures}).
The partial derivatives of the recourse solution with respect to the reserves $\rp$ and $\rm$ are identically zero.
Partial derivatives of the solution to the generation levels $\piw$ with respect to the variables $\po$ and $\alpha$ are computed by differentiating Eqns.~\eqref{eqn:targets},~\eqref{eqn:final_bal}, and~\eqref{eqn:smooth_sat} and solving the resulting linear system of sensitivities.
Partial derivatives of the phase angle solution $\tiw$ with respect to $\po$ and $\alpha$ are computed by differentiating Eqns.~\eqref{eqn:flowbal} and solving the resulting linear system.
We summarize these relationships in Appendix \ref{app:partialderiv}. 

\subsection{Solving the smooth approximation using PSG}

We use the projected stochastic gradient (PSG) method of~\citep{ghadimi2016mini,davis2018stochastic} to obtain stationary solutions to the smooth approximation\footnote{An alternative is to use sample average approximation (SAA) for approximating the solution to the smooth approximation, which can also exploit its decomposable structure}.
We first note that assumption (A2) of~\citep{davis2018stochastic} holds since we assume that the conditions of Theorem~\ref{thm:consis_uniq} hold.
Furthermore, the smooth saturation function~$g_{\tau}$ 
is continuously differentiable with Lipschitz continuous gradient.
The sensitivities of the recourse solutions with respect to the first-stage variables are also Lipschitz continuous.
Hence, the objective function of our smooth approximation is continuously differentiable with Lipschitz continuous gradient.
Because the first-stage feasible region is compact, assumption $\overline{\text(A3)}$ of~\citep{davis2018stochastic} also holds and the PSG method is guaranteed to converge to stationary solutions. A detailed description of the algorithm is provided in Appendix~\ref{app:stochapprox}.

\section{Alternative Models}
\label{sec:alternatives}

We compare the solution of the smooth approximation with the solutions from two alternative models that determine candidate first-stage decisions using the affine policy in Eqn.~\eqref{eqn:affine_policy} instead of the saturation model in Eqns.~\eqref{eqn:targets} and~\eqref{eqn:smooth_sat}. 
\subsubsection{Conservative Affine Policy (CAP) Model}
The first model we compare against is inspired by~\citep{bienstock2014chance,zhang2016distributionally,
vrakopoulou2017optimal}. This model enforces individual generator limits using chance constraints with maximum violation allowances $\varepsilon_{gen}$, while still using penalty terms to limit line violations. The motivation behind this model is avoid the need to consider saturation effects by ensuring that the generation amounts from the affine policy very rarely hit the generator limits. 
\subsubsection{Generator Penalty (GP) Model} The second alternative we consider does not directly include a constraint on the violation probability, but rather penalizes the expected violation of the generator limits by the generation levels determined by the affine policy (cf.~\citep{roald2016optimal}) using the terms $\gamma_{gen}\max\{0,\piw - \pmaxiw,\pminiw - \piw\}^2$, $i \in \G$, in the recourse objective for a penalty coefficient $\gamma_{gen} > 0$. 
In our computational experiments, we investigate whether it is possible to choose $\gamma_{gen}$ such that the GP model yields good solutions to the true Problem~\eqref{eqn:stochopf}.

Both of these alternative models are solved using a nonlinear solver to solve a sample average approximation.
Appendix~\ref{app:formulations} provides detailed descriptions of these models. 

\section{Computational experiments}
\label{sec:computexp}

\subsection{Modeling and implementation details}

We set $\pminiw \equiv 0$, $\forall i \in \G$, $\pol = 0$, $\pou_i = \pmax_i$, $\forall i \in \Rg$, and $\pou_i = \expect{\pmaxiw} + 5\stdev{\pmaxiw}$, $\forall i \in \W$.
For the reserve bounds, we set $\rpu = \rmu = \pou$.
We use $\Gres = \G$ by default, and let $\varepsilon = \min\left\lbrace 0.001,\frac{0.01}{\abs{\G}} \right\rbrace$. For each generator $i \in \G$, we use smoothing parameter $\tau_{sat} = 10^{-4} \left( \pmaxiw - \pminiw \right)$. We assume that the generation limits are wide enough for relatively complete recourse to hold.
The cost functions are specified as $f_{1,i}(z) = c_i z$, $f_{2,i}(z) = f_{3,i}(z) = c_i c_{res} z$, $\forall i \in \Rg$, for input unit generation costs $c_i$ and reserves cost factor $c_{res} = 1.5$.
For the wind generators, we use $f_{1,i}(z) \equiv 0$ (marginal) generation costs and reserve cost functions $f_{2,i}(z) = f_{3,i}(z) = \left(\min_{k \in \Rg} c_k\right) c_{wind} c_{res} z$, $\forall i \in \W$, where $c_{wind} =0.1$ is the relative cost factor for wind reserves.
The penalty functions are specified as $q_{1,i}(z) \equiv 0$, $q_{2,i}(z) = \gamma_{res} f_{2,i}(g^+_{\tau_{pos}}(z))$, and $q_{3,i}(z) = \gamma_{res} f_{3,i}(-g^-_{\tau_{pos}}(z))$, $\forall i \in \G$, 
where $\gamma_{res} = 10$ is the penalty for reserves beyond the scheduled capacity, $\tau_{pos} = 10^{-4}$ is a smoothing parameter, $g^+_{\tau_{pos}}$ is the smooth approximation to the $(\cdot)_+$ function defined by:
{
\small
\[ g^+_{\tau_{pos}}(z) := \tau_{pos} \log\Bigl( 1 + \exp\bigl( \dfrac{z}{\tau_{pos}} \bigr) \Bigr),\]
}%
and $g^-_{\tau_{pos}}(z) := -g^+_{\tau_{pos}}(-z)$ is the smooth approximation to $(z)_{-}$.
Finally, $q_{4,ij}(z) := \gamma_{line}\max\{0,\abs{z} - \delta_{ij}\fmax_{ij}\}^2$, where $\fmax_{ij}$ is the $(i,j)^{\text{th}}$ line flow limit, $\delta_{ij} \equiv 0.95$, and $\gamma_{line} > 0$ is the line flow penalty coefficient that is varied.


We solve $500$ scenario SAAs of the comparison models CAP and GP (which can be expressed as convex quadratic programs) to determine candidate first-stage solutions.
We use a solution from the GP model with $\gamma_{gen} = 20$ as the initial guess $x_1$ for our smooth approximation model.
The quality of the solutions obtained using all approaches are evaluated on the true model~\eqref{eqn:stochopf} (i.e., including reserve saturation) using a common independent Monte Carlo sample of size~$10^5$.

The code and data of the test instances are available at \url{https://github.com/rohitkannan/DCOPF-reserve-saturation}.
Our codes are written in Julia~0.6.2~\citep{bezanson2017julia}, and use Gurobi~7.5.2~\citep{gurobi} to solve convex programs through the JuMP~0.18.2 interface~\citep{dunning2017jump}. We use IPOPT~3.12.8~\citep{wachter2006implementation} in situations where Gurobi encounters numerical difficulties.
All computational tests were conducted on a Surface Book 2 laptop running Windows 10 Pro with a $1.90$~GHz four core Intel~i7 CPU, $16$~GB of RAM.



\begin{figure}
    \centering
     \begin{subfigure}[b]{0.48\textwidth}
         \centering
         \includegraphics[width=1.0\textwidth]{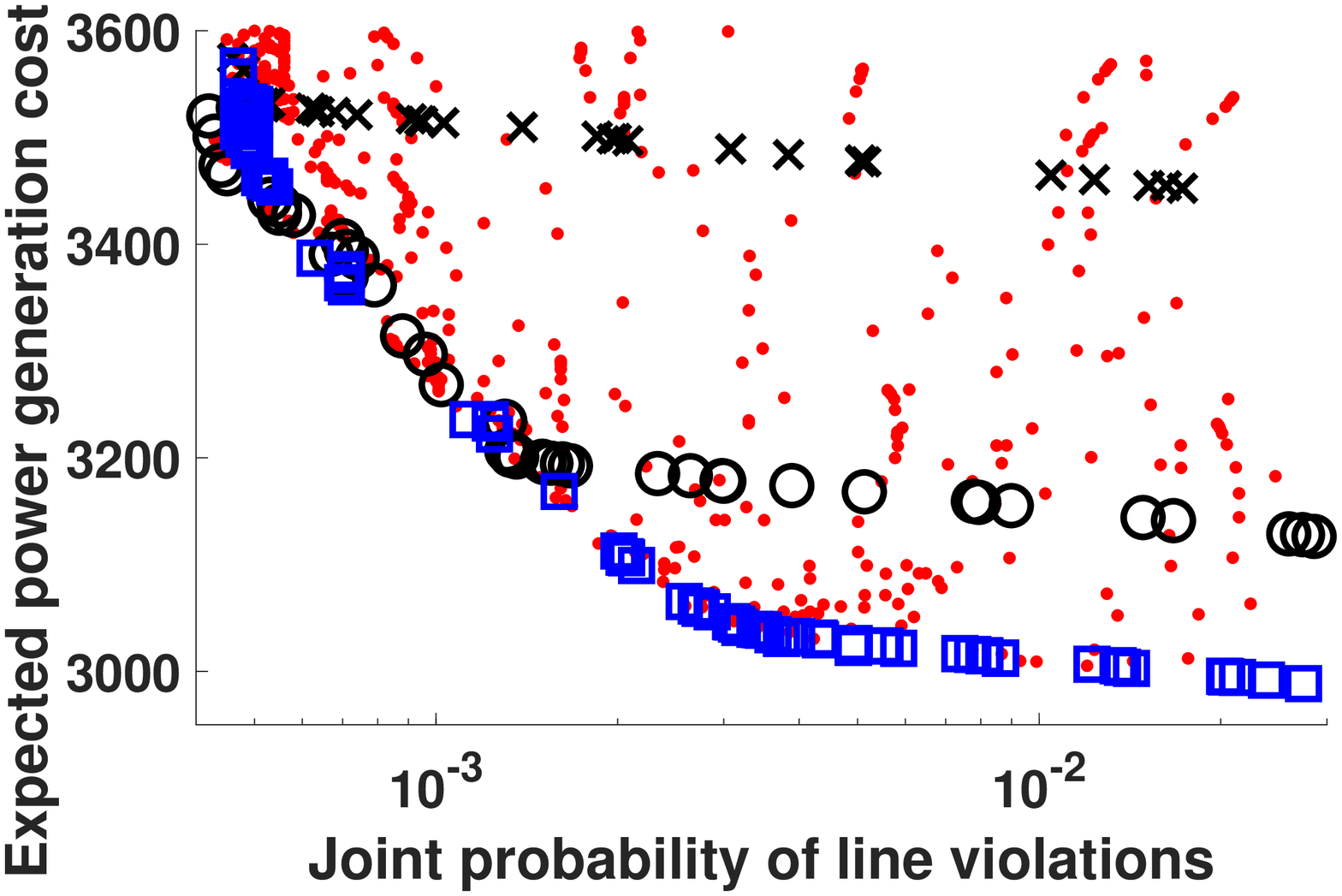}
     \end{subfigure}
     \begin{subfigure}[b]{0.48\textwidth}
         \centering
         \includegraphics[width=1.0\textwidth]{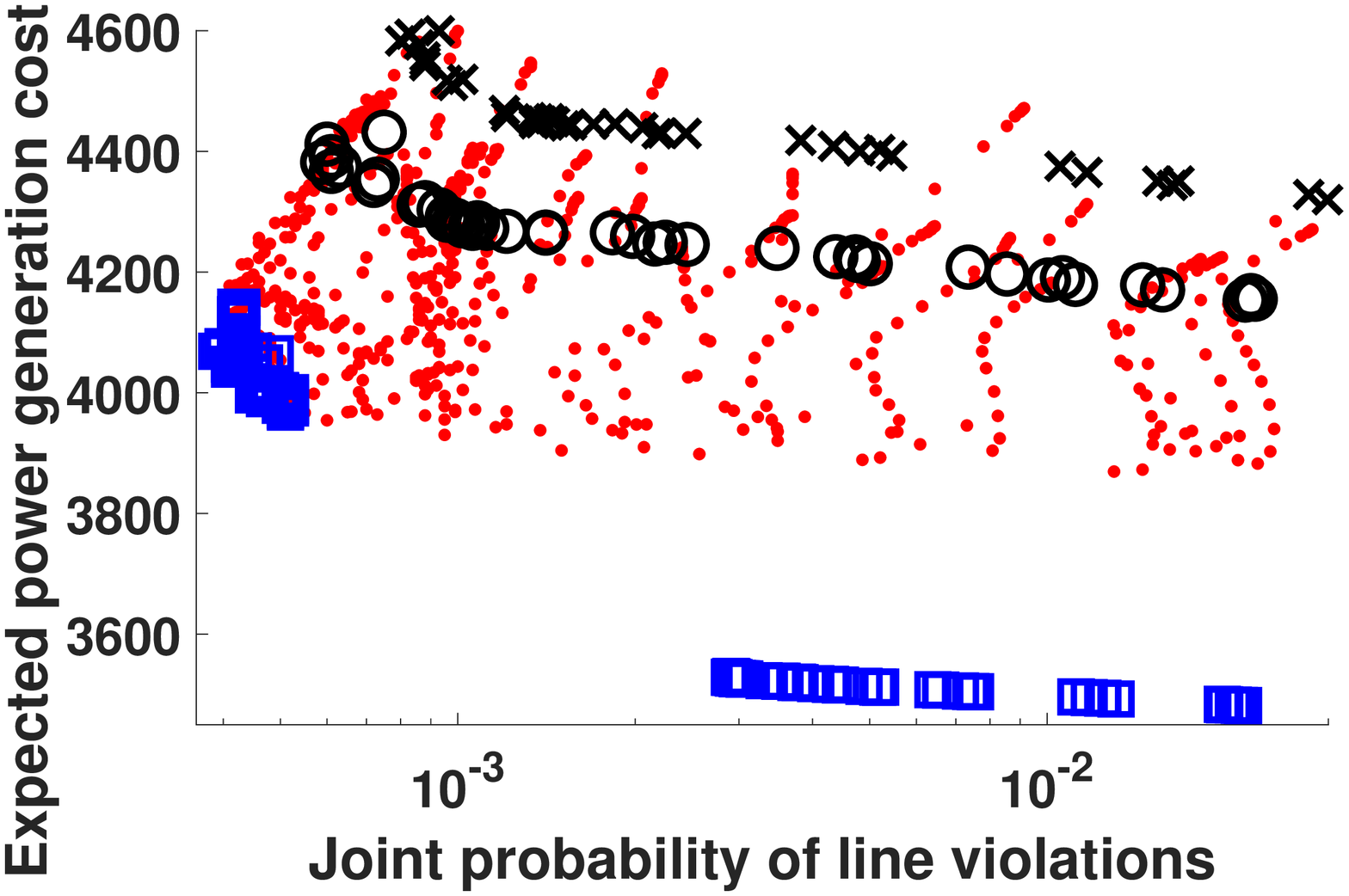}
     \end{subfigure}    
     \begin{subfigure}[b]{0.48\textwidth}
         \centering
         \includegraphics[width=1.0\textwidth]{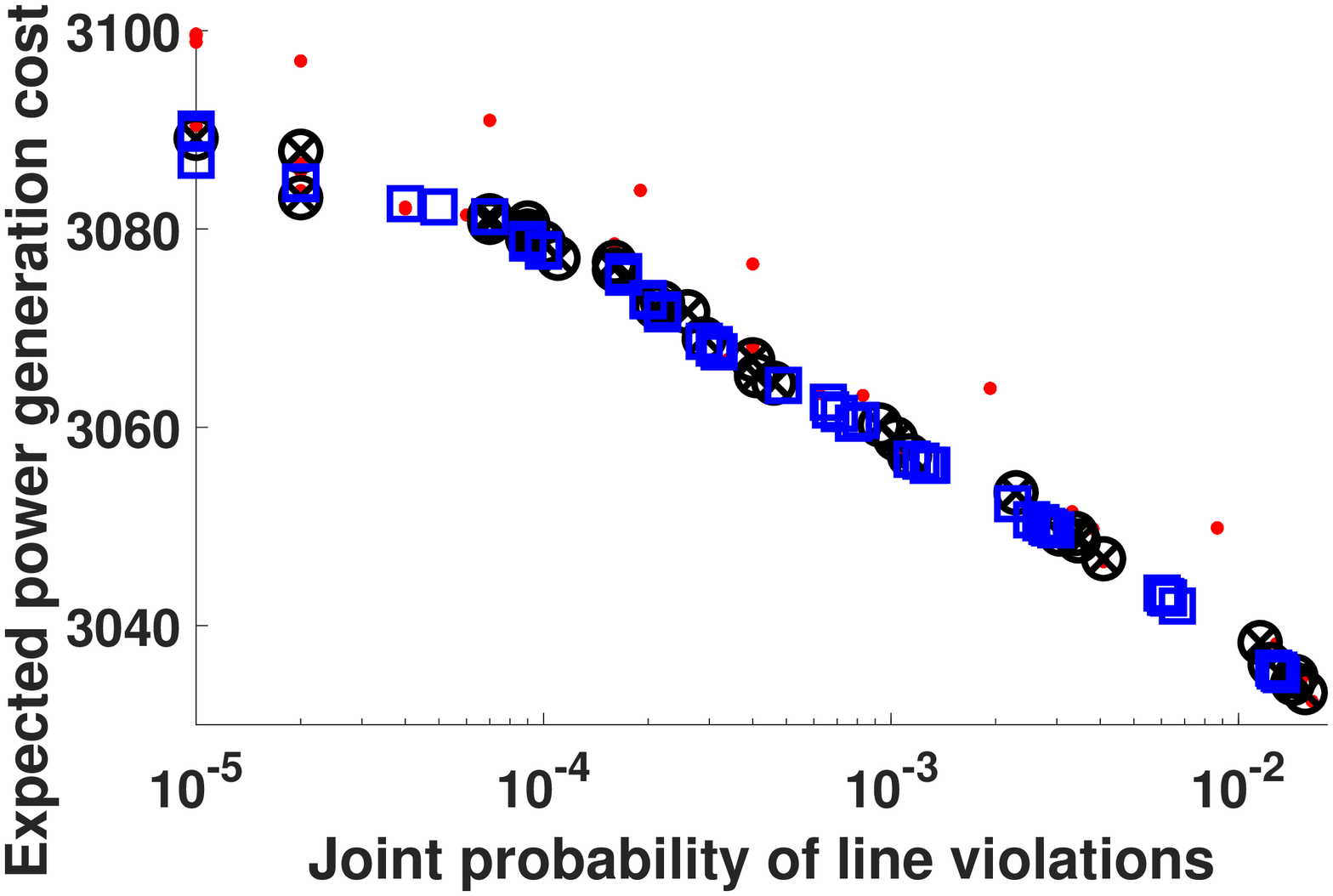}
     \end{subfigure}      
     \caption{\small (Top to bottom) Pareto plots for the 6-bus system generated using five replicates for \casenum{1}, \casenum{2}, and \casenum{3}. \textbf{\color{blue}Blue squares:} Smooth approximation \:\: \textbf{\color{red}Red dots:} GP solution \:\: \textbf{Black crosses:} CAP $10^{-3}$ \:\: \textbf{Black circles:} CAP $10^{-2}$}
    \label{fig:pareto_6bus} 
\end{figure}

\begin{figure*}
    \centering
     \begin{subfigure}[b]{0.48\textwidth}
         \centering
         \includegraphics[width=1.0\textwidth]{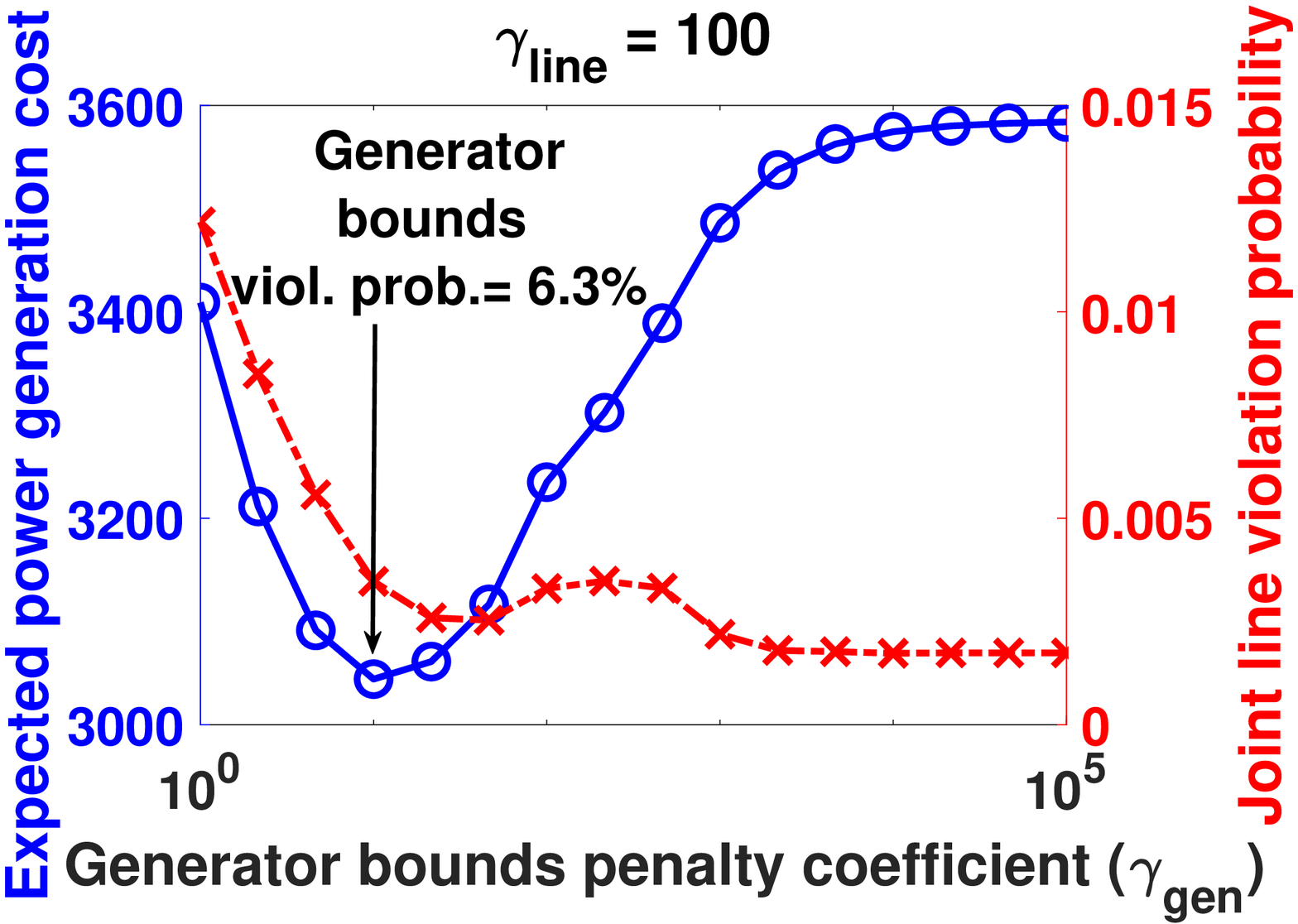}
     \end{subfigure}
     \hfill
     \begin{subfigure}[b]{0.48\textwidth}
         \centering
         \includegraphics[width=1.0\textwidth]{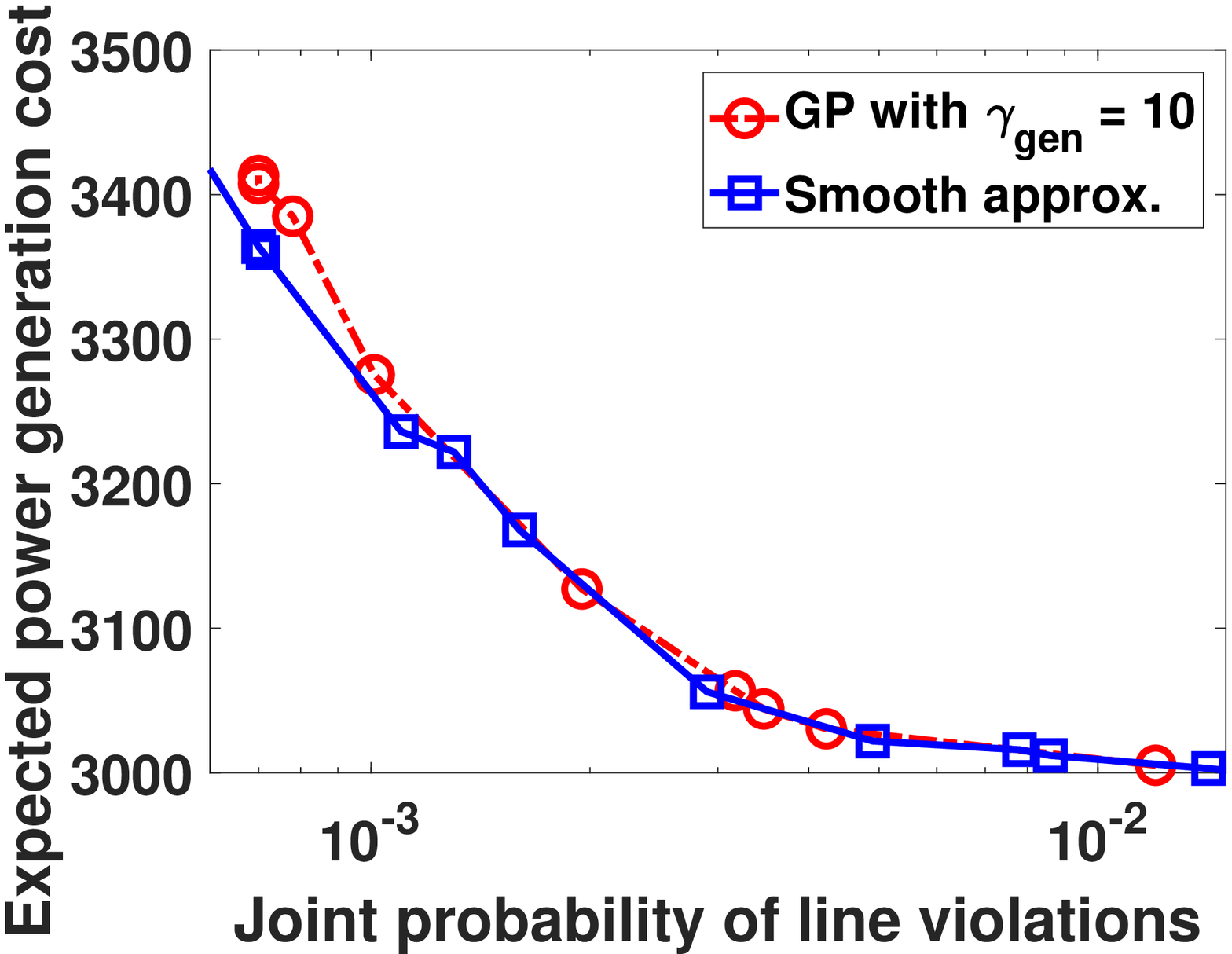}
     \end{subfigure}    
     \caption{\small (Left) Expected cost of power generation (\textbf{\color{blue} blue circles}) and joint line violation probability (\textbf{\color{red} red crosses}) versus generator penalty coefficient $\gamma_{gen}$ for \casenum{1} of the 6-bus system generated with fixed line penalty coefficient $\gamma_{line} = 100$. (Right) Expected cost of power generation versus joint line violation probability for \casenum{1} of the 6-bus system with fixed generator penalty coefficient $\gamma_{gen} = 10$ and varying line penalty coefficient $\gamma_{line}$.}
    \label{fig:gp_6bus} 
    \vspace*{-0.15in}
\end{figure*}

\begin{table}
\caption{\small Comparison of the lowest cost solutions with joint line violation probability $\leq 0.5\%$ for the 6-bus model \casenum{1}.}
\label{tab:soln_comp_windres}
\resizebox{\columnwidth}{!}{%
\begin{tabular}{c|c|cc|c|c}
        & Expected   & \multicolumn{2}{|c|}{First-stage cost}   & Gen.    & Wind \\
Model   & total cost & Gen. & Res.            & viol.                & util. \% \\ 
\hline
SA
& 3049  & 2713  & 115  & -         & 70   \\
CAP $10^{-3}$ 
& 3505  & 3188  & 226  & $<$0.002 & 45.2  \\
CAP $10^{-2}$ 
& 3186  & 2991  & 126  & $<$0.01& 55.9  \\
GP        
& 3043  & 2689  & 110  & 0.069     & 71.4    \\
\hline
\end{tabular}
}
\end{table}

\begin{table}
\caption{\small Comparison of the lowest cost solutions with joint line violation probability $\leq 0.5\%$ for the 6-bus model \casenum{2}.}
\label{tab:soln_comp_nowindres}
\resizebox{\columnwidth}{!}{%
\begin{tabular}{c|c|cc|c|c}
        & Expected   & \multicolumn{2}{|c|}{First-stage cost}   & Gen.    & Wind \\
Model   & total cost & Gen. & Res.            & viol.                & util. \% \\ 
\hline
SA
& 3546  & 1514  & 1675  & -         & 100   \\
CAP $10^{-3}$ 
& 4446  & 2797  & 1127  & $<$0.001 & 69.1  \\
CAP $10^{-2}$ 
& 4268  & 2614  & 1122  & $<$0.01 & 76.7  \\
GP        
& 3904  & 2133  & 1186  & 0.567     & 96.8    \\
\hline
\end{tabular}
}
\end{table}%

\subsection{Case Study I: 6-bus system}
Our 6-bus example (with $\abs{\G} = 3$) is based on~\url{http://motor.ece.iit.edu/data/6bus_Data_ES.pdf}.
We recast `generator G2' as a wind generator, and consider normally distributed loads and wind generator capacities with average wind output equal to half the average load.
We consider three cases: \\
\casenum{1}: wind generators \textit{can} provide reserves,\\ \casenum{2}: wind generators \textit{do not} provide reserves (but are allowed to spill wind without cost, with $\alpha_{wind} = 0.1\varepsilon$), and \\
\casenum{3}: wind generators are non-dispatchable (i.e., they act like negative loads).\\
We assume that the standard deviation of the wind output and the loads are $10\%$ of the average for the first two cases, but only $5\%$ of the average for the third case to ensure relatively complete recourse. For these three cases, we compare the solution obtained with the smooth approximation, as well as the solutions from the CAP and GP models.

\subsubsection{Pareto plots for the three algorithms} 
To provide a complete picture of the performance of solutions that can be obtained from each algorithm, we present Pareto plots to display the quality of solutions obtained across different parameter values. To generate the Pareto plot, we do a parameter sweep for the tuning parameters of each model. The line violation penalty parameters are changed between $\gamma_{line} = \text{logspace}(1,5,17)$ for our smooth approximation (except for \casenum{2}, where we use $21$ values of $\gamma_{line}$ between $10^{1}$ and $10^{5}$), and between $\gamma_{line} = \text{logspace}(1,5,9)$ for the GP and CAP models. For the GP model we also do a parameter sweep for the generator violation penalties with $\gamma_{gen} = \text{logspace}(0,5,16)$. For the CAP model, we consider two different violation probabilities for generator chance constraints, viz., $\varepsilon_{gen} = [10^{-3}, 10^{-2}]$. Since the solutions depend on the samples and are therefore random, we create five replications for each parameter combination. 
For each solution, we calculate the expected cost of power generation (including cost of reserves and reserve penalties) and the joint probability that \textit{any} line flow limit is violated by evaluating the system behavior based on the true policy (which includes reserve saturation) on an independent sample. Note that we do not include any assessment of the generator violation probability, since this probability is zero in the true model.



Fig.~\ref{fig:pareto_6bus}
shows the Pareto plots for the three different algorithms and the three different cases, with expected generation cost plotted against the expected joint violation probability for the line flows. We plot solutions for the smooth approximation (blue squares), the GP model (red dots) and the CAP model with two different values for $\varepsilon_{gen}$ (black circles and crosses).

We observe from the Pareto plots that our smooth approximation always provides nearly non-dominated solutions for all three cases. 
The solutions obtained with the GP model are not concentrated along the Pareto front, as generation violation penalties that are either too large or too small lead to higher-than-necessary cost. Solutions from the CAP model provide a different Pareto front with a higher cost than the smooth approximations, though the solutions coincide with the smooth approximation for smaller values of the violation probability.

Beyond these general behaviors, the models compare differently between different cases. In \casenum{1}, careful parameter tuning allows the GP model to find points along the Pareto curve, see Fig.~\ref{fig:gp_6bus}. In \casenum{2}, there is a gap between the lowest cost solutions obtained with the GP model and the Pareto curve found with the smooth approximation. The smooth approximation is able to find lower cost solutions when the probability of line violations is not set too low. 
The lack of data points in the Pareto curve for the smooth approximation model in \casenum{2} between line violation probabilities of $6 \times 10^{-4}$ and $2 \times 10^{-3}$ is a result of using a linear weighting approach for balancing the expected generation and violation costs for the nonconvex Problem~\eqref{eqn:stochopf}, see Chapter~3 of~\citep{ehrgott2006multicriteria}.
In \casenum{3}, the solutions of all three algorithms cluster along the Pareto front. 

\subsubsection{Detailed comparison of differences} To analyze the cause of these differences, we investigate some of the solutions in more detail. For each algorithm and cases $1$ and $2$, we list the results for the lowest cost solution with a joint line violation probability $\leq0.5\%$ in Tables~\ref{tab:soln_comp_windres} and~\ref{tab:soln_comp_nowindres}. We list the total expected cost, the first stage scheduled generation and reserve capacity costs, the out-of-sample joint violation probabilities of the lines,
and the joint violation probability of the generators if we would have considered an affine control policy (only calculated for the CAP and GP models). We also list the expected utilization of wind energy, as a percentage of total available wind power. 

 In \casenum{1} we observe that the solutions from GP and the smooth approximation have very similar total costs and utilization of the wind energy. The smooth approximation invests more in both generation and reserve capacity in the first stage, which is balanced by paying lower penalties in the second stage. Interestingly, the optimal choice of tuning parameters for the GP solution leads to a relatively high violation probability for the generators at $6.9\%$ (also see Fig.~\ref{fig:gp_6bus}). If we enforce a lower violation probability, as is done in the CAP model and has typically been done in literature (see e.g. \cite{bienstock2014chance}), the total expected cost increases significantly and the utilization of wind energy drops. 

For \casenum{2}, where wind generators are not allowed to provide reserves, we observe that the total expected cost is significantly lower for the smooth approximation than for both the GP and CAP models. The smooth approximation has a lower first-stage generation cost (indicating high dispatch levels for the cheap wind power), but invests more in procuring reserves (that can make up for overestimates in the wind generation). This leads to full utilization of the available wind power. In comparison, the GP model schedules less wind power in the first stage, leading to a higher cost and lower wind utilization. Interestingly, the affine policy in the best GP solution violates the generator limits with more than $50\%$ probability. This also explains why the CAP solutions, where the generation violation probability is explicitly limited, leads to much higher total expected cost (and much lower wind power utilization) than the other two models. 

Finally, we do not compare the solutions in \casenum{3} as they are very similar for all three algorithms. In this case, the solutions balance the cost of scheduling more power from the less expensive generator with paying penalties for violating the line constraints. The reserve activation is happening at the more expensive generator, which is far away from saturation. The generators never violate their limits even with an affine control policy and the models are therefore the same.

\begin{figure}
    \centering
     \begin{subfigure}[b]{0.48\textwidth}
         \centering
         \includegraphics[width=1.0\textwidth]{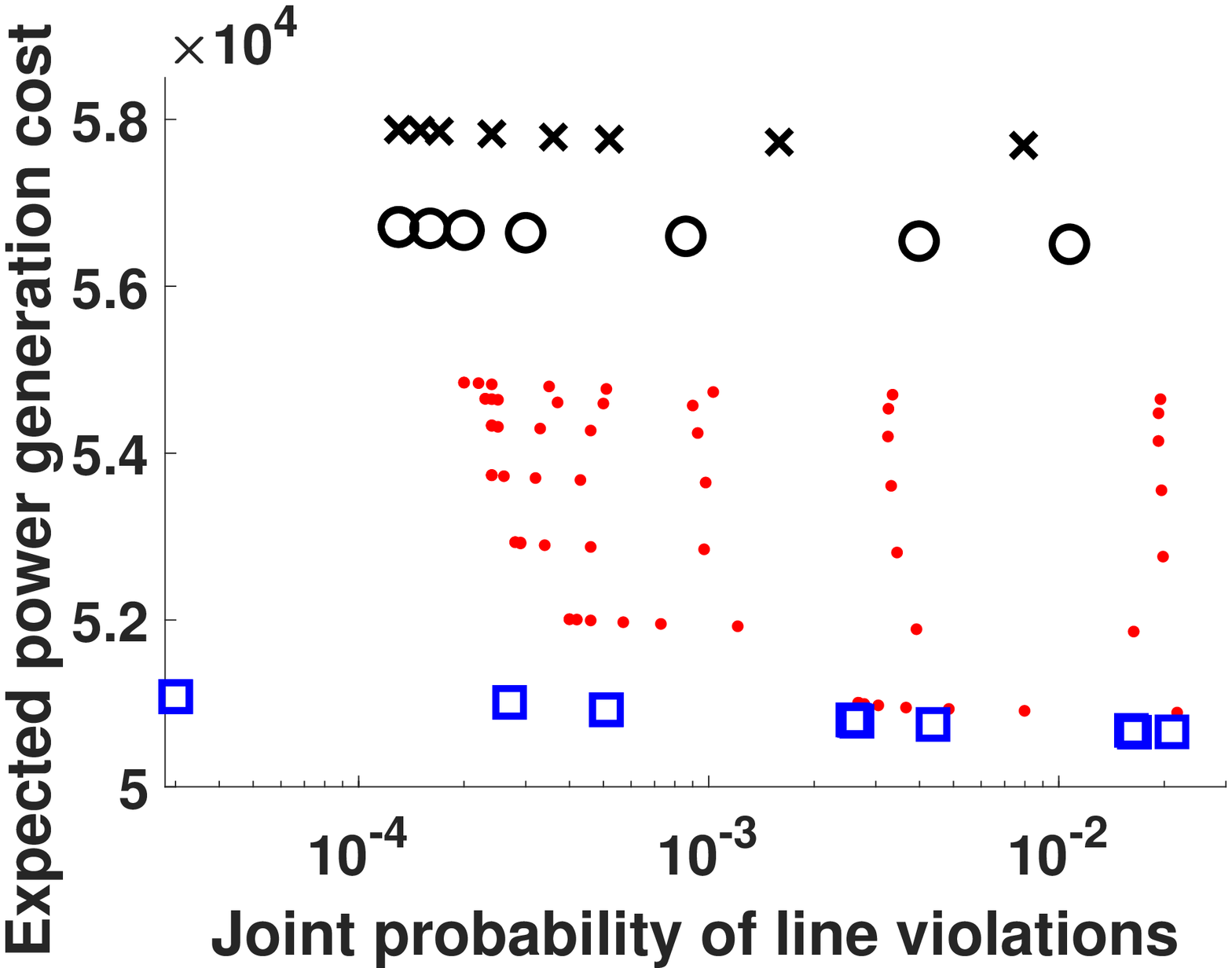}
     \end{subfigure}
     \begin{subfigure}[b]{0.48\textwidth}
         \centering
         \includegraphics[width=1.0\textwidth]{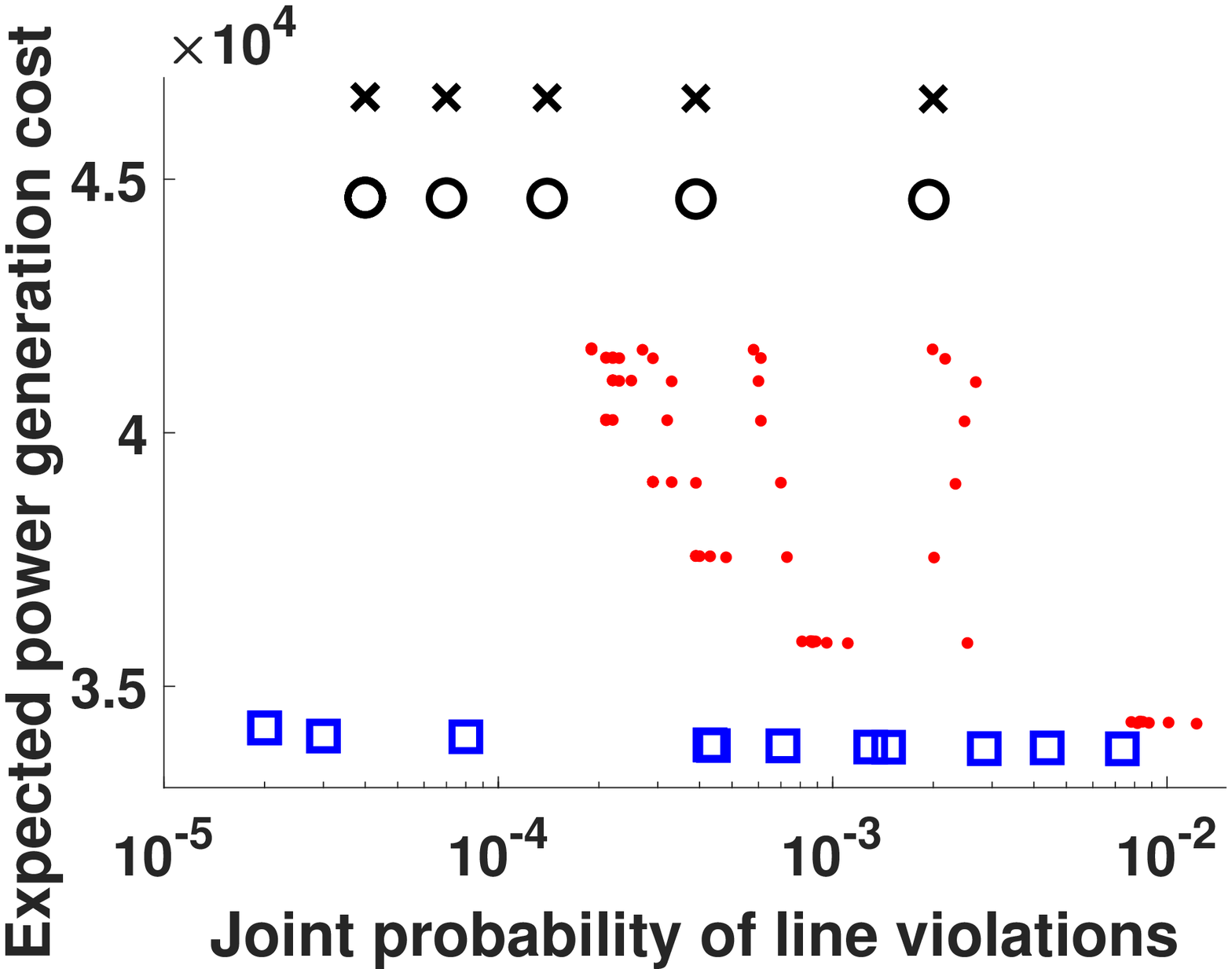}
     \end{subfigure}    
     \begin{subfigure}[b]{0.48\textwidth}
         \centering
         \includegraphics[width=1.0\textwidth]{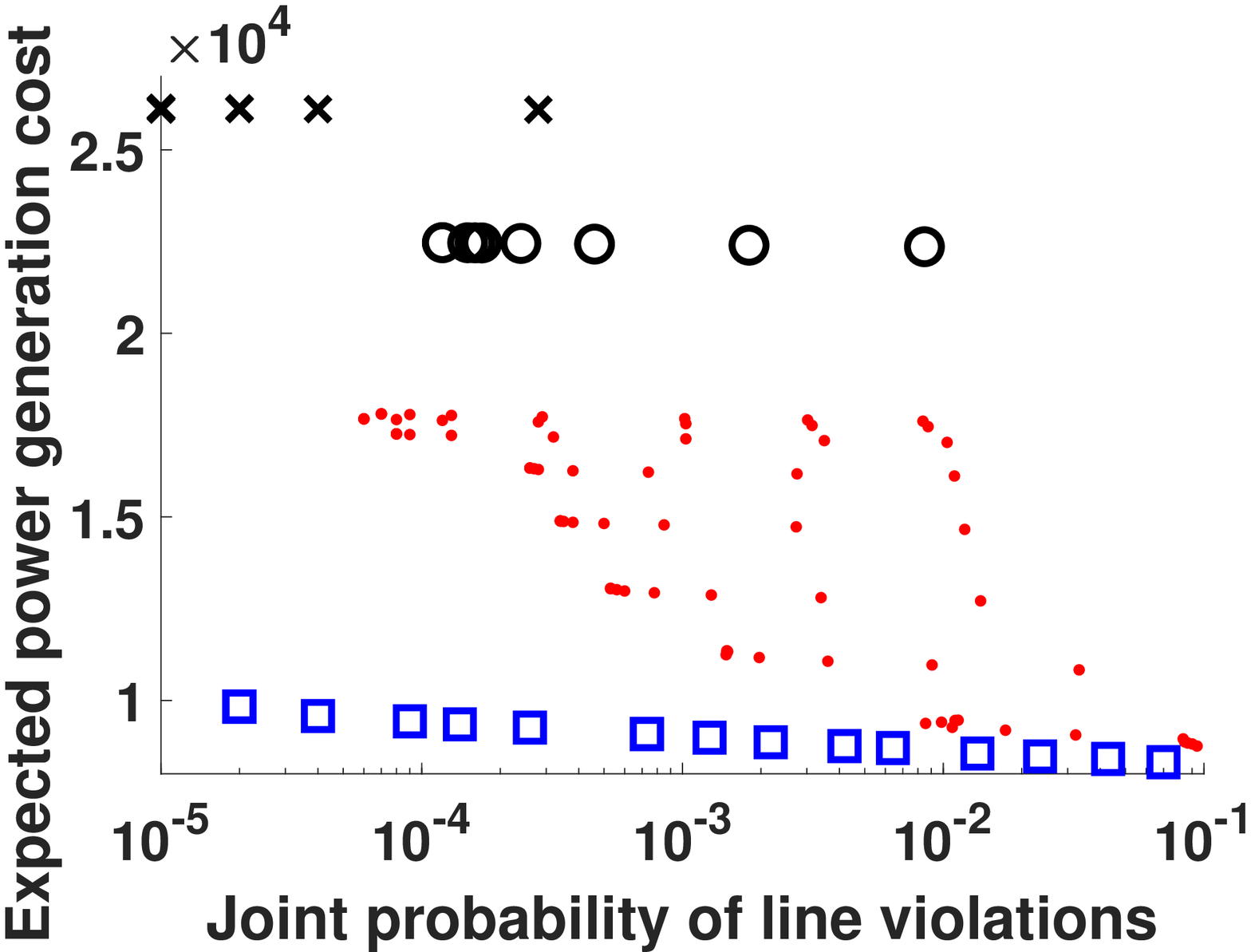}
     \end{subfigure}      
     \caption{\small (Top to bottom) Pareto plots for the 118-bus system generated using five replicates for $25\%$, $50\%$, and $100\%$ wind penetration. \textbf{\color{blue}Blue squares:} Smooth approximation \: \textbf{\color{red}Red dots:} GP solution \: \textbf{Black crosses:} CAP $10^{-5}$ \: \textbf{Black circles:} CAP~$10^{-4}$}
    \label{fig:pareto_118bus} 
\end{figure}

\begin{figure}
    \centering
     \begin{subfigure}[b]{0.48\textwidth}
         \centering
         \includegraphics[width=1.0\textwidth]{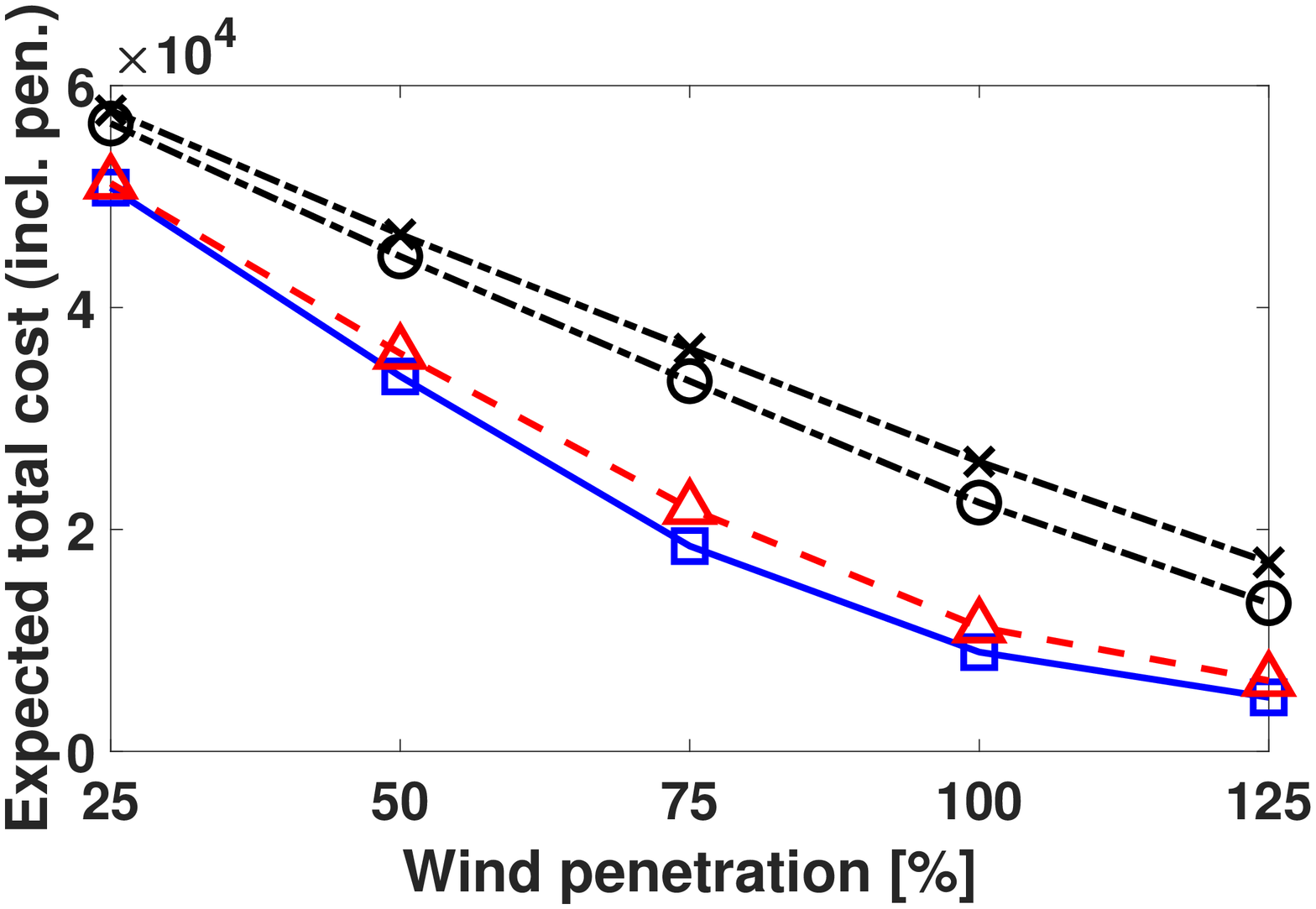}
     \end{subfigure}
     \begin{subfigure}[b]{0.48\textwidth}
         \centering
         \includegraphics[width=1.0\textwidth]{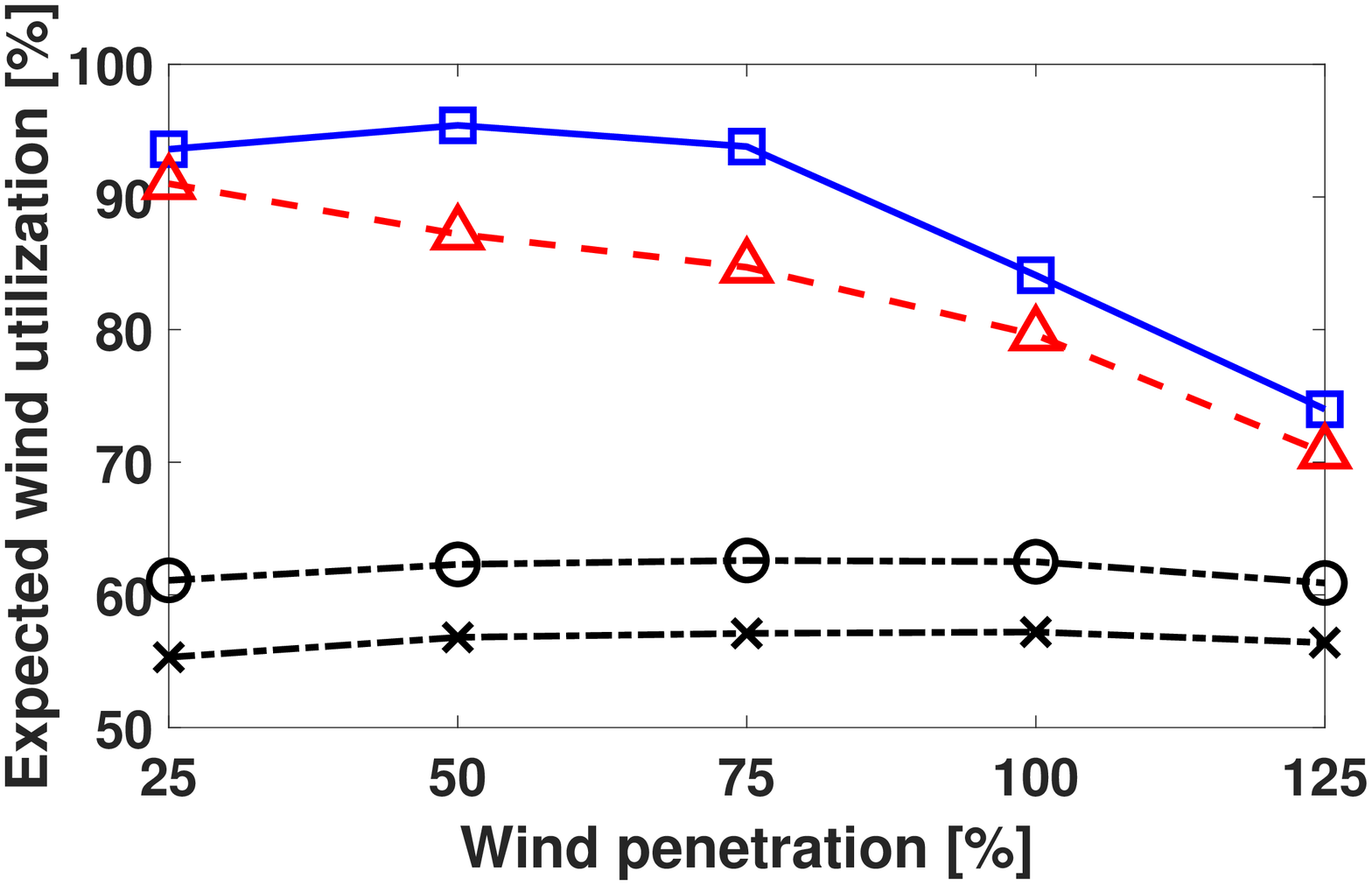}
     \end{subfigure}    
     \begin{subfigure}[b]{0.48\textwidth}
         \centering
         \includegraphics[width=1.0\textwidth]{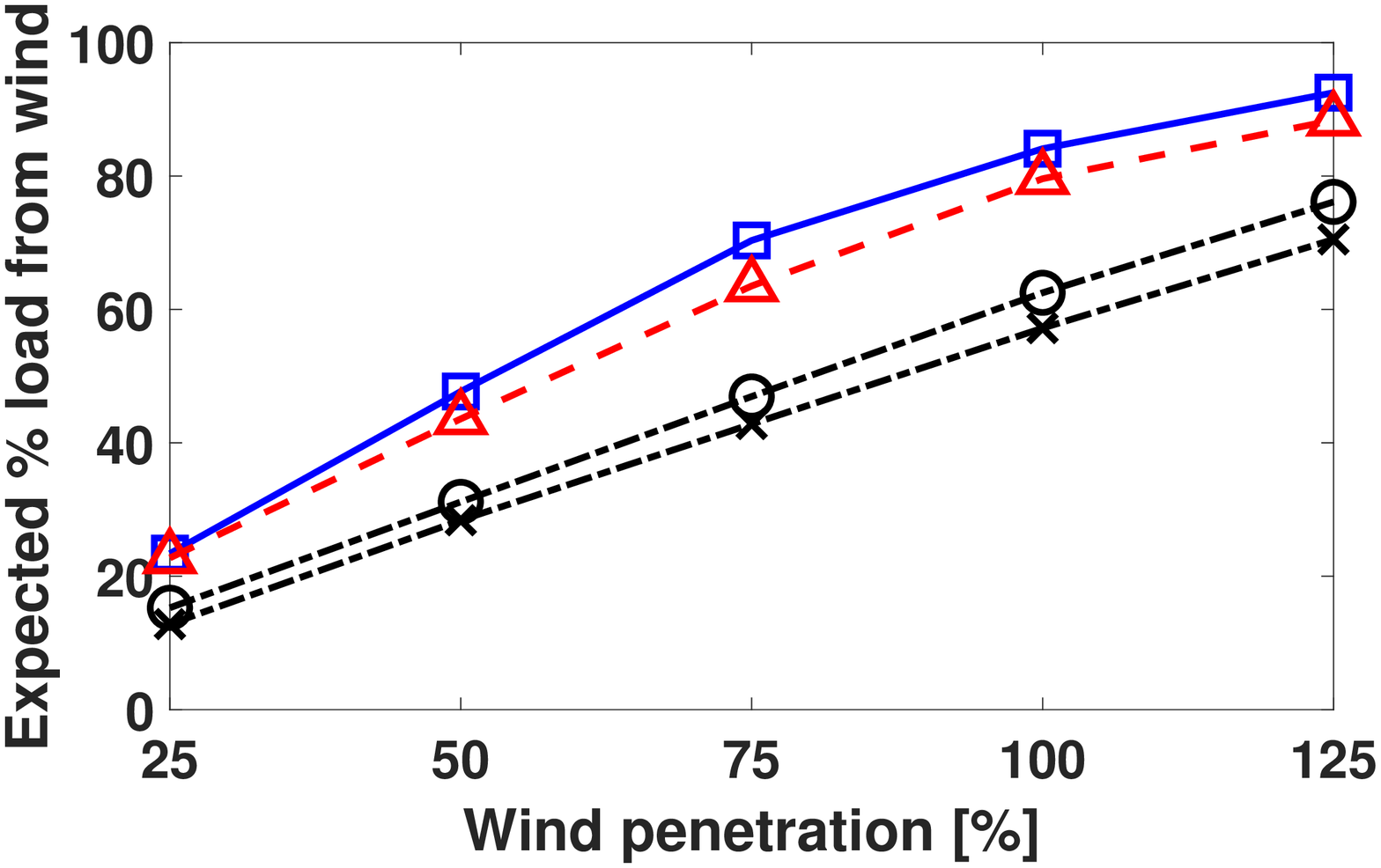}
     \end{subfigure} 
     \caption{\small Summary of solutions for the different models in the 118-bus case with varying wind penetration levels. \textbf{\color{blue}Blue squares: }Smooth approximation~~~\textbf{\color{red}Red triangles:} GP solution~~~\textbf{Black crosses:} CAP $10^{-5}$~~~\textbf{Black circles:} CAP $10^{-4}$}
    \label{fig:windpenetration} 
    \vspace*{-0.15in}
\end{figure}

\subsection{Case Study II: 118-bus system}

In the second part of our case study, we consider the more realistic test case based on the IEEE 118-bus system from \url{http://motor.ece.iit.edu/data/JEAS_IEEE118.doc} with modifications suggested in~\citep{roald2016optimal}, including the addition of $25$ wind generators to the $54$ regular generators, increasing the average demand by $50\%$, and reducing the line flow limits by $25\%$.
We again consider normally distributed loads and wind generator capacities, and consider five different levels of wind penetration: average wind output $=$ $25\%$, $50\%$, $75\%$, $100\%$, or $125\%$ of the average system load.
To obtain appropriate parameter values for the algorithms, we run a similar parameter sweep as for the 6-bus test case, but with $\gamma_{gen} = \text{logspace}(0,4,9)$ and $\varepsilon_{gen} = [10^{-5},10^{-4}]$ to limit computational effort. 
Fig.~\ref{fig:pareto_118bus}
shows the Pareto plots for the three different algorithms and the three different cases for wind penetration levels of $25\%$, $50\%$, and $100\%$, with expected generation cost plotted against the expected joint violation probability for the line flows. We plot solutions for the smooth approximation (blue squares), the GP model (red dots) and the CAP model with two different values for $\varepsilon_{gen}$ (black circles and crosses).
We compare the solution with lowest expected total cost and joint line flow violation probability $\leq 0.5\%$ for each algorithm and each wind level penetration. 
The expected total cost, expected wind utilization and expected fraction of total system load served by wind power is calculated using a Monte Carlo simulation with $10^5$ samples and the results are plotted in Fig.~\ref{fig:windpenetration}. 
It takes roughly $1.5$ minutes on average to solve the GP and CAP models and roughly $7$ minutes on average to solve the smooth approximation model and generate a point on the Pareto curve for this case.

The smooth approximation outperforms the other methods in all aspects. The cost is lower, the wind utilization is higher and more load is served by the wind generators. The solution obtained by a properly tuned affine GP model achieves similar performance as the smooth approximation (though it has a $30\%$ higher expected cost in some cases), while both of the CAP methods perform significantly worse.

\section{Conclusion and future work}
\label{sec:concl}

We propose a stochastic DC optimal power flow model with reserve saturation. Specifically, our model assumes that generators follow an affine control policy until they reach a generation limit, at which point they operate at that limit. The model is a two-stage stochastic program with nonconvex, nonsmooth second stage constraints, and we propose a stochastic approximation method to solve a smooth approximation. 
We empirically observe that our model yields solutions that outperform those obtained from a model that constrains the affine control policy to rarely violate generation limits. On the other hand, using a model that penalizes expected violation of generator limits can sometimes yield competitive solutions with a well-tuned choice of the penalty parameter.

Extensions to our model that would be interesting to investigate in future work include constraining the probability or expected violation of line limits rather than penalizing violation of line limits in the objective (which will provide additional data points to fill the gap in the Pareto curve for the smooth approximation model in \casenum{2}), and using an AC power flow model in place of the DC model.



{
\footnotesize
\section*{References}
\begingroup
\renewcommand{\section}[2]{}%
\bibliographystyle{unsrtnat}
\bibliography{main}
\endgroup
}

\appendices

\section{Proof of Theorem~\ref{thm:consis_uniq}}
\label{app:proof}
If $\sigmadw \not\in D_F(\po,\alpha,\random)$, then Eqns.~\eqref{eqn:targets},~\eqref{eqn:saturation}, and~\eqref{eqn:final_bal} are inconsistent since Eqns.~\eqref{eqn:targets} and~\eqref{eqn:saturation} together imply
{
\begin{align*}
&\:\: \sum_{i \in \Resa} \pminiw + \sum_{j \not\in \Resa} \pbojw - \sum_{k \in \D} d_k \\
\leq &\:\: \sum_{i \in \G} \piw - \sum_{k \in \D} d_k \\
\leq &\:\: \sum_{i \in \Resa} \pmaxiw + \sum_{j \not\in \Resa} \pbojw - \sum_{k \in \D} d_k,
\end{align*}
}
which makes the satisfaction of Eqn.~\eqref{eqn:final_bal} impossible.

We now show that the system of equations~\eqref{eqn:targets},~\eqref{eqn:saturation}, and~\eqref{eqn:final_bal} has a solution whenever $\sigmadw \in D_F(\po,\alpha,\random)$.
For a chosen value of the slack reserves~$\sw$, denote the value of $\piw$ obtained using Eqns.~\eqref{eqn:targets} and~\eqref{eqn:saturation} by $\hat{p}_i(s)$ (we omit dependence on the first-stage variables and $\random$ for simplicity). 
Note that $\hat{p}_i(s)$ is a monotonically nondecreasing continuous function of $s$ for each $i \in \G$, which implies that $\sum_{i \in \G} \hat{p}_i(s)$ is a monotonically nondecreasing continuous function of $s$. 
Furthermore, we have from Eqns.~\eqref{eqn:targets} and~\eqref{eqn:saturation} that generators with a nonzero participation factor will eventually hit their bounds for small/large enough chosen values of~$s$, i.e., there exists $M > 0$ large enough for which $\hat{p}_i(-M) = \pminiw$ and $\hat{p}_i(M) = \pmaxiw$, $\forall i \in \Resa$, which implies $\sum_{i \in \G} \hat{p}_i(-M) = \sum_{i \in \Resa} \pminiw + \sum_{j \not\in \Resa} \pbojw$ and $\sum_{i \in \G} \hat{p}_i(M) = \sum_{i \in \Resa} \pmaxiw + \sum_{j \not\in \Resa} \pbojw$. 
From the intermediate value theorem applied to $\sum_{i \in \G} \hat{p}_i(\cdot)$, there exists $\hat{s} \in [-M,M]$ such that $\sum_{i \in \G} \hat{p}_i(\hat{s}) = \sum_{j \in \D} \left[ d_j + \dtjw\right]$ for any $\sigmadw \in D_F(\po,\alpha,\random)$. 
Therefore, the system of equations~\eqref{eqn:targets},~\eqref{eqn:saturation}, and~\eqref{eqn:final_bal} has a solution for the variables $\left(\ptiw, \piw, \sw \right)$ whenever $\sigmadw \in D_F(\po,\alpha,\random)$.

When $\sigmadw \in \text{int}(D_F(\po,\alpha,\random))$, we have from Eqn.~\eqref{eqn:final_bal} that there exists a generator that has not hit its generation limits, i.e., $\exists j \in \Resa$ such that $\pminjw < p_j(\random) < \pmaxjw$ at a solution to Eqns.~\eqref{eqn:targets},~\eqref{eqn:saturation}, and~\eqref{eqn:final_bal}. This implies that the sum $\sum_{i \in \G} \hat{p}_i(\cdot)$ is monotonically (strictly) increasing in a neighborhood of $\sw$ around the above solution, which establishes its uniqueness since $\sum_{i \in \G} \hat{p}_i(s)$ is a monotonically nondecreasing continuous function of $s$. 
The argument for the `only if' part is similar.
\qed

\vspace*{-0.05in}

\section{Calculation of Partial Derivatives}
\label{app:partialderiv}

The linear system below is solved to obtain the partial derivatives of the recourse solution with respect to the first-stage decision variable $q$, where $q$ is a placeholder for either $\po_l$ or $\alpha_l$, $l \in \G$:
\begin{align}
\frac{\partial \pt_i}{\partial q}(\random) &= \frac{\partial \po_i}{\partial q} + \left(\sw + \sigmadw\right)\frac{\partial \alpha_i}{\partial q} + \alpha_i \frac{\partial s}{\partial q}(\random), \label{eq:partial1} \\
\frac{\partial p_i}{\partial q}(\random) &= \frac{\partial g_{\tau_{sat}}}{\partial \pt_i}(\ptiw;\pminiw,\pmaxiw) \frac{\partial \pt_i}{\partial q}(\random), \\
\sum_{k \in \G} \frac{\partial p_k}{\partial q}&(\random) = 0, \quad \frac{\partial \theta_1}{\partial q}(\random) = 0
\end{align}
\begin{align}
\sum_{j \: : \: (i,j) \in \E} &\beta_{ij} \left[ \frac{\partial \theta_i}{\partial q}(\random) - \frac{\partial \theta_j}{\partial q}(\random) \right] = \frac{\partial p_i}{\partial q}(\random). \label{eq:partial4}
\end{align}

\section{Stochastic Approximation Algorithm Details}
\label{app:stochapprox}

Algorithm~\ref{alg:stochgrad} presents a basic version of our PSG workflow. In this algorithm, we use $x := (\po,\rp,\rm,\alpha)$ as the set of all first-stage variables, and define $X$ to be the set of $x$ that satisfy Eqns.~\eqref{eqn:init_bal} to~\eqref{eqn:participation}. The operator $\mathrm{Proj}_X(y)$ returns the point in $X$ that has smallest Euclidean distance to $y$. For simplicity the algorithm is described with a fixed step length $\gamma$, but we use a variation of AdaGrad~\citep{duchi2011adaptive} for determining step lengths.

\begin{algorithm}
\caption{PSG algorithm for solving the smooth approx.}
\label{alg:stochgrad}
{
\begin{algorithmic}[1]

\State \textbf{Input}: Initial guess $x_1 \in X$, number of iterations $T \in \mathbb{N}$, mini-batch size $K \in \mathbb{N}$, and step length $\gamma > 0$.

\For{$t = 1, \cdots, T$}
\For{$k=1,\cdots,K$}
\parState{Let $\random_k$ be a random observation of $\omega$.}
\parState{Solve Eqns.~\eqref{eqn:targets},~\eqref{eqn:flowbal}, and~\eqref{eqn:smooth_sat} to obtain $\swk$ and $p^{\text{T}}(\random_k)$ for the given $x_t$ and $\random_k$.}
\parState{Solve Eqns.~\eqref{eq:partial1}-\eqref{eq:partial4} with $\random_k$, $\swk$, $p^{\text{T}}(\random_k)$ and $x_t$ and use the chain rule to get a stochastic gradient $\hat{g}_k$ of the objective of~\eqref{eqn:stochopf} at iterate $x_{t}$.}
\EndFor
\State Let $x_{t+1} = \mathrm{Proj}_X \left(x_{t} - \gamma \frac{1}{K}\sum_{k=1}^K\hat{g}_k\right)$.
\parState{Estimate the objective of Problem~\eqref{eqn:stochopf} using an independent sample of~$\random$, and check termination criteria.}
\EndFor

\State \textbf{Output}: Iterate with the smallest estimated objective.

\end{algorithmic}
}
\end{algorithm}

\section{Details of the Formulations}
\label{app:formulations}

We explicitly write out the three formulations considered for the computational experiments below (parameter settings are listed in Sec.~\ref{sec:computexp}).

\noindent\textbf{Smooth approximation}:
\begin{alignat*}{2}
&\uset{\rm, \alpha}{\min_{\po,\rp,}} && \displaystyle\sum_{i \in \G} \left[ c_i \po_i + \bar{c}_i \left(\rp_i + \rm_i \right)\right] + Q_1(\po, \rp, \rm, \alpha) \nonumber \\
&\quad \text{s.t.} && \:\: \text{Constraints}~\eqref{eqn:init_bal} \text{ to}~\eqref{eqn:participation},
\end{alignat*}
where $\bar{c}_i = c_i c_{res}$, for $i \in \Rg$, and $\bar{c}_i =
c_{wind} c_{res} (\min_{k \in \Rg}{c_k})$, for  $i \in \W$, $c_i = 0$, $\forall i \in \W$, $Q_1(\po,\rp,\rm,\alpha) = \expectation{\random}{q_1(\po,\rp,\rm,\alpha,\random)}$ denotes the expected second-stage costs with $q_1(\po,\rp,\rm,\alpha,\random) :=$
{
\small
\begin{alignat*}{1}
\uset{\sw,\tw}{\min_{\pw,\ptw,}} \nonumber &\sum_{i \in \G} \gamma_{res} \bar{c}_i g^+_{\tau_{pos}}(\piw - \po_i - \rp_i) + \\
&  \sum_{i \in \G} \gamma_{res} \bar{c}_i g^+_{\tau_{pos}}(\po_i - \piw - \rm_i ) + \nonumber \\
&\hspace*{-0.05in}\sum_{(i,j) \in \E} \max\left\lbrace 0, \abs*{\beta_{ij} \left[ \tiw - \tjw \right]} - \delta_{ij} \fmax_{ij} \right\rbrace^2 \\
\text{s.t.} \quad & \text{Constraints}~\eqref{eqn:targets},~\eqref{eqn:flowbal}, \text{ and}~\eqref{eqn:smooth_sat}. \nonumber
\end{alignat*}
}

\noindent\textbf{Conservative affine policy (CAP) model}
{
\small
\begin{alignat*}{2}
&\uset{\rm, \alpha}{\min_{\po,\rp,}} && \displaystyle\sum_{i \in \G} \left[ c_i \po_i + \bar{c}_i \left(\rp_i + \rm_i \right)\right] + Q_2(\po, \rp, \rm, \alpha) \nonumber \\
&\quad \text{s.t.} && \:\: \text{Constraints}~\eqref{eqn:init_bal} \text{ to}~\eqref{eqn:participation}, \\
& && \:\: \prob{\po_i + \alpha_i \Sigma_d(\random) \geq \pminiw} \geq 1 - \varepsilon_{gen}, \:\: \forall i \in \G, \\
& && \:\: \prob{\po_i + \alpha_i \Sigma_d(\random) \leq \pmaxiw} \geq 1 - \varepsilon_{gen}, \:\: \forall i \in \G,
\end{alignat*}
}
where $Q_2(\po,\rp,\rm,\alpha) = \expectation{\random}{q_2(\po,\rp,\rm,\alpha,\random)}$ denotes the expected second-stage costs with $q_2(\po,\rp,\rm,\alpha,\random) :=$
{
\small
\begin{alignat*}{1}
\uset{\pw,\tw}{\min} \nonumber &\sum_{i \in \G} \gamma_{res} \bar{c}_i \left( \left(\piw - \po_i - \rp_i \right)_{+} + \right. \\
&\quad\quad\quad\quad\quad\quad \left. \left(\po_i - \piw - \rm_i \right)_{+} \right) + \nonumber \\
&\sum_{(i,j) \in \E} \max\left\lbrace 0, \abs*{\beta_{ij} \left[ \tiw - \tjw \right]} - \delta_{ij} \fmax_{ij} \right\rbrace^2 \\
\text{s.t.} \quad & \text{Constraints}~\eqref{eqn:affine_policy} \text{ and}~\eqref{eqn:flowbal}.
\end{alignat*}
}

\noindent\textbf{Generator penalty (GP) model}
\begin{alignat*}{2}
&\uset{\rm, \alpha}{\min_{\po,\rp,}} && \displaystyle\sum_{i \in \G} \left[ c_i \po_i + \bar{c}_i \left(\rp_i + \rm_i \right)\right] + Q_3(\po, \rp, \rm, \alpha) \nonumber \\
&\quad \text{s.t.} && \:\: \text{Constraints}~\eqref{eqn:init_bal} \text{ to}~\eqref{eqn:participation},
\end{alignat*}
where $Q_3(\po,\rp,\rm,\alpha) = \expectation{\random}{q_3(\po,\rp,\rm,\alpha,\random)}$ denotes the expected second-stage costs with $q_3(\po,\rp,\rm,\alpha,\random) :=$
{
\footnotesize
\begin{alignat*}{1}
\uset{\pw,\tw}{\min} \nonumber &\sum_{i \in \G} \gamma_{res} \bar{c}_i \left( \left(\piw - \po_i - \rp_i \right)_{+} + \right. \\
&\quad\quad\quad\quad\quad\quad \left.\left(\po_i - \piw - \rm_i \right)_{+} \right) + \nonumber \\
&\sum_{(i,j) \in \E} \max\left\lbrace 0, \abs*{\beta_{ij} \left[ \tiw - \tjw \right]} - \delta_{ij} \fmax_{ij} \right\rbrace^2 + \\
&\sum_{i \in \G} \gamma_{gen}\max\left\lbrace 0,\piw - \pmaxiw,\pminiw - \piw \right\rbrace^2 \\
\text{s.t.} \quad & \text{Constraints}~\eqref{eqn:affine_policy} \text{ and}~\eqref{eqn:flowbal}.
\end{alignat*}
}

\end{document}